%% file: ms.tex
\documentclass[a4paper,10pt]{amsart}
\usepackage[T1]{fontenc}
\usepackage[utf8]{inputenc}
\usepackage{amsmath}
\usepackage{amsthm}
\usepackage{amssymb}
\usepackage{mathtools}
\usepackage{booktabs}
\usepackage{longtable}
\usepackage{pdflscape}
\usepackage{color}

\theoremstyle{plain}
\newtheorem{theorem}{Theorem}
\newtheorem{lemma}{Lemma}

\theoremstyle{remark}

\newtheorem{example}{Example}
\newlength{\temp}
\title{Local Type II metrics with holonomy in $\mathrm{G}_2^*$}
\author{Christian Volkhausen}
\address{Christian Volkhausen, Institut für Mathematik und Informatik, Universität Greifs\-wald, Walther-Rathenau-Straße 47, 17489 Greifswald, Germany}
\email{christian.volkhausen@uni-greifswald.de}
\urladdr{\scriptsize https://math-inf.uni-greifswald.de/en/department/about-us/employees/christian-volkhausen/}
\thanks{I am grateful to Ines Kath for her useful comments and advice on this article.}
\begin{document}
\begin{abstract}
A list of possible holonomy groups with indecomposable holonomy representation contained in the exceptional, non-compact Lie group $\mathrm{G}_2^{*}$ was provided by Fino and Kath. The classification is due to the corresponding holonomy algebras and divided into Type I, II and III, depending on the dimension of the socle being 1, 2 or 3, respectively. It was also shown by Fino and Kath that all algebras of Type I, and by the author that all of Type III are indeed realizable as holonomy algebras by metrics with signature (4,3). This article proves that this is also true for all Type II algebras. Thus, there exists a realization by a metric for all indecomposable holonomy groups contained in $\mathrm{G}_2^{*}$.
\end{abstract}%
\maketitle
\section{Introduction}\label{Sec:Introduction}
The holonomy of locally nonsymmetric semi-Riemannian manifolds is described by Berger's list \cite{Berger1955} if the holonomy representation, i.e., the natural representation of the holonomy group on the tangent space, is irreducible. Here, irreducibility means that there is no invariant proper subspace of the holonomy representation. Additionally, Cartan classified Riemannian symmetric spaces by their holonomy groups \cite{Cartan1926}. Thus, holonomy groups of Riemannian manifolds are completely classified.

In the case of pseudo-Riemannian manifolds the situation is more complicate. Since the holonomy representation may have isotropic subspaces, it is not sufficient to classify only the irreducible representations. Instead, it is necessary to consider the much larger class of indecomposable ones. Indecomposability here means that the holonomy representation does not leave invariant any proper nondegenerate subspace. 

While holonomy groups of pseudo-Riemannian manifolds with irreducible representation are also known by Berger's list, there are only a few results in the classification of indecomposable ones. Leistner provided a complete classification for Lorentzian manifolds \cite{Leistner2007}. Galaev classified holonomy groups of Lorentz-Kähler, i.e., those with Kähler metrics of index 2, and Einstein pseudo-Riemannian manifolds \cite{Galaev2018b,Galaev2018a}. Aside from these classifications some partial results are known (cf. \cite{Bezvitnaya2011,Fino2015,Graham2012}).

In this article we deal with pseudo-Riemannian manifolds of signature $(4,3)$. By Berger's list, the holonomy group of such a manifold is either generic $SO(4,3)$ or $\mathrm{G}_2^{*}$ if the holonomy representation is irreducible. Therefore, it is natural to ask which proper subgroups $H$ of the exceptional, non-compact Lie group $\mathrm{G}_2^{*}$ are holonomy groups with indecomposable but not irreducible representation. The corresponding Lie algebras $\mathfrak{h} \subset \mathfrak{g}_2^*$ are necessarily Berger algebras, i.e., fulfill Berger's first criterion. A complete classification of these subalgebras up to conjugation in $SO(4,3)$ is provided by Fino and Kath  \cite{Fino2016}.  Additionally, Fino and Kath distinguish those subalgebras by the dimension of their maximal, semisimple subrepresentation on the tangent space, the so called socle. An algebra is called of Type I, II or III if the socle has dimension 1, 2 or 3, respectively. While these results provide a classification of  Berger subalgebras in $\mathfrak{g}_2^{*}$, nothing is said about their realization as a holonomy algebra. It is potentially a highly nontrivial task to prove that a given Berger algebra is realizable as a holonomy algebra by some metric. For instance, it took about thirty years to answer this question for the groups $\mathrm{G}^{c}_2$ and $Spin(7)$ on Berger's list  \cite{Bryant1987}. 

In the case of proper subalgebras of $\mathfrak{g}_2^{*}$, Fino and Kath were able to realize all Type I algebras by appropriate metrics \cite{Fino2017} and the author did the same for Type III \cite{Volkhausen2019}. Here, we want to answer the same question for Type II in an affirmative way. Thus, we prove  
\begin{theorem}\label{Thm:MainTheorem}
Each indecomposable Berger algebra $\mathfrak{h}\subset \mathfrak{g}_2^{*}$ of Type II is realizable as a holonomy algebra by some metric with signature $(4,3)$.
\end{theorem}
Together with earlier results on Type I and III this also proves the general theorem
\begin{theorem}
Each indecomposable Berger algebra $\mathfrak{h}\subset \mathfrak{g}_2^{*}$ is realizable as a holonomy algebra by some metric with signature $(4,3)$.
\end{theorem}
In section 2 we summarize some properties of the group $\mathrm{G}_2^{*}$ and holonomy algebras of Type II. We keep this section brief, since most of the algebraic background is exhaustively discussed. We recommend \cite{Fino2016,Fino2017,Kath1998,Volkhausen2019} for details. Therefore, we concentrate on the differences between Type I and III on the one hand and Type II on the other. The proof of Theorem \ref{Thm:MainTheorem} itself is constructive. The construction approach is covered by section 3. It is fairly the same as used in case of Type I and III. Finally, a metric is explicitly indicated in the appendix for each algebra we have to consider. 
\section{Type II holonomy algebras}\label{Sec:2}
In this section we provide a brief review of Type II algebras contained in $\mathfrak{g}_2^{*}$, the Lie algebra of $\mathrm{G}_2^{*}$. The original description is due to Fino and Kath \cite{Fino2016}, whose notation is adopted with minor additions. Finally, we state the classification theorem for Type II algebras.

The Lie algebra $\mathfrak{g}_2^{*}$ of the exceptional, real, non-compact Lie group $\mathrm{G}_2^{*}$, as well as the group itself, is well known from classical Lie theory. Classically, this group is characterized as the stabilizer of a real, generic 3-form $\omega \in \Lambda^{3}V$. Here $V$ is a 7-dimensional vector space and `generic' means that $\omega$ has open $GL(7,V)$-orbit. This characterization traces back to Engel, who used a similar approach in order to characterize the complex group $\mathrm{G}_2$ \cite{Engel1900}.

There are some further characterizations equivalent to the classical one. We can consider $\mathrm{G}_2^{*}$ as the stabilizer of a non-isotropic spinor in the real spinor representation of $Spin(4,3)$. Also, $\mathrm{G}_2^{*}$ stabilizes a certain cross-product on $\mathbb{R}^{4,3}$. Starting with one of these characterizations we can obtain the others by suitable calculations in a straightforward way. For details we refer to \cite{Kath1998}.

Now, let $M$ be a smooth, connected, simply connected manifold. We assume simply connectedness, since we are interested in reduced holonomy groups only. Denote the tangent space $T_xM$ at some point $x \in M$ with $V$ and identify $V\cong \mathbb{R}^{7}$ by choosing the basis $(b_i)$ of $V$ as
\begin{align*}
b_i = e_i \;\mathrm{for}\; i\in\{1,\ldots,4\},\quad  b_5 = e_7, \quad b_6 = e_5, \quad b_7 = e_6 \hspace{.5em},
\end{align*}
where $e_i$ denote the canonical basis of $\mathbb{R}^{7}$. We use this specific basis mainly for the sake of consistency with \cite{Fino2016}. Note, the basis differs from that used for Type I and III in \cite{Fino2017,Volkhausen2019}, which causes the main differences in notation among the articles. Now, with respect to this basis we describe $\mathrm{G}_2^{*}$ as the stabilizer of the generic 3-form $\omega \in \Lambda^{3}V^{*}$ given by
\begin{align*}
\omega = \sqrt{2} \left( -b^{157} + b^{236}\right) - b^4 \wedge \left(b^{16}-b^{27}-b^{35}\right)  \hspace{.5em},
\end{align*}
where $b^i$ denotes the dual basis and $b^{ijk} = b^i \wedge b^j \wedge b^k$. By  Bryant \cite{Bryant1987}, $\omega$ induces an inner product via the relation
\begin{align*}
 \langle X,Y\rangle &= \frac{1}{6} (\iota_{X}(\omega) \wedge \iota_{Y}(\omega) \wedge \omega) \in \Lambda^{7}V \cong \Lambda^{7} \mathbb{R}^{7} \cong \mathbb{R}\hspace{.5em}.
\end{align*}
for any vectors $X,Y \in V$. In terms of the chosen basis of $V$ this can be expressed as
\begin{align}\label{Eq:InnerProduct}
\langle \cdot,\cdot\rangle = 2  (b^1 \cdot b^6 + b^2 \cdot b^7 + b^3 \cdot b^5) - b^4 \cdot b^4 \hspace{.5em}.
\end{align}
With respect to the inner product we regard $V \cong \mathbb{R}^{4,3}$ and consider $\mathrm{G}_2^{*}$ as a subgroup of $SO(4,3)$.

We now turn to the Lie algebra $\mathfrak{g}_2^{*}$ of $\mathrm{G}_2^{*}$. With respect to the basis $b_i$, elements of $\mathfrak{g}_2^{*}$ are given by matrices of the form
\begin{align}\label{Eq:g2_elements}
\begin{pmatrix}
s_1+s_4&-s_{10}&s_9&\sqrt 2s_6&-s_{12}&0&-s_{11}\\
-s_8&s_1&s_2&\sqrt 2s_9&s_6&s_{11}&0\\
s_7&s_3&s_4&\sqrt 2s_{10}&0&s_{12}&-s_6\\
 \sqrt 2s_5&\sqrt 2
s_7&\sqrt 2s_8&0&\sqrt 2s_{10}&\sqrt 2s_6&\sqrt 2s_9\\ 
-s_{14}&s_5&0&\sqrt 2s_8&-s_4&-s_9&-s_2\\
0&s_{13}&s_{14}&\sqrt 2s_5&-s_7&-s_1-s_4&s_8\\
 -s_{13}&0&-s_5&\sqrt 2s_7&-s_3&s_{10}&-s_1
\end{pmatrix} \hspace{.5em},
\end{align}
where $s_i \in \mathbb{R}$. 
Furthermore, $\mathfrak{g}_2^{*}$ is a $|2|$-graded Lie algebra, i.e., as a vector space it decomposes into a direct sum
\begin{align*}
\mathfrak{g}_{2}^{*} =\mathfrak{g}_{-2} \oplus \mathfrak{g}_{-1} \oplus \mathfrak{g}_0 \oplus \mathfrak{g}_{1} \oplus \mathfrak{g}_{2}
\end{align*} 
such that $[\mathfrak{g}_i,\mathfrak{g}_j] \subset \mathfrak{g}_{i+j}$. Using the representation \eqref{Eq:g2_elements} and the root diagram of $\mathfrak{g}_2^{*}$ known from classical Lie theory (cf. \cite{Fulton2004}), subspaces $\mathfrak{g}_i$ are easily determined by commuting several elements of $\mathfrak{g}_2^{*}$ with each other. One finds that $\mathfrak{g}_0$ is spanned by a Cartan subalgebra, the shorter primitive root vector and its negative. We also find $\dim \mathfrak{g}_{\pm 1} = 4$ and $\dim \mathfrak{g}_{\pm 2} = 1$. Since the longer primitive root vector is contained in $\mathfrak{p}_2 = \mathfrak{g}_0 \oplus \mathfrak{g}_1 \oplus \mathfrak{g}_2$, the subalgebra $\mathfrak{p}_2$ contains a Borel subalgebra, hence is parabolic. Indeed, it is maximal parabolic and of dimension 9. With respect to \eqref{Eq:g2_elements} $\mathfrak{p}_2$ is given by 
\begin{align*}
\mathfrak{p}_2 = \left\lbrace m \in \mathfrak{g}_2^{*}\, \middle| \,s_i = 0 \;\mathrm{for}\; i = 3,5,7,13,14 \right\rbrace 
\end{align*} 
One can show that the corresponding maximal parabolic subgroup $P_2 \subset \mathrm{G}_2^{*}$ is the stabilizer of certain 2-plains in $\mathbb{R}^{4,3}$. Note, the algebras of Type I and III are contained in the other maximal parabolic subalgebra $\mathfrak{p}_1 \subset \mathfrak{g}_2^{*}$. The corresponding maximal parabolic subgroup $P_1 \subset \mathrm{G}_2^{*}$ stabilizes an isotropic line in $\mathbb{R}^{4,3}$. Thus, aside from notational issues there are algebraic differences between Type I and III algebras on the one hand and Type II on the other. We will see that holonomy algebras of Type II are contained in $\mathfrak{p}_2$.

Let $z=(z_1,\ldots,z_4)\in \mathbb{R}^4$, $A = \begin{pmatrix} a_1 & a_2 \\ a_3 & a_4 \end{pmatrix}\in \mathfrak{gl}(2,\mathbb{R})$, $c \in \mathbb{R}$. Then, define
\begin{align*}
\sigma(z) &= \begin{pmatrix} z_2 & \sqrt{2} z_3 & z_4 \\ z_1 & \sqrt{2} z_2 & z_3 \\ \end{pmatrix}\;, \qquad
\sigma(z)^* = \begin{pmatrix} -z_4 & -z_3 \\ \sqrt{2}z_3 & \sqrt{2} z_2 \\ -z_2 & -z_1 \end{pmatrix} \\
\rho(A) &= \begin{pmatrix}a_1-a_4 & -\sqrt{2} a_2 & 0 \\ -\sqrt{2}a_3 & 0 & -\sqrt{2} a_2 \\ 0 & -\sqrt{2}a_3 & -a_1+a_4\end{pmatrix}\;, \quad U(c) = \begin{pmatrix} 0 & -c \\ c & 0 \end{pmatrix} \hspace{.5em},
\end{align*}
which provides a suitable notational frame for the description of Type II algebras as given in the following lemma by Fino and Kath
\begin{lemma}[Type II holonomy algebras, \cite{Fino2016}]\label{Lem:TypeIIAlgebras}
 If $\mathfrak{h}$ is of Type II, then there exists a basis $b_1,\ldots,b_7$ of $V$ such that 
\begin{align*}
\omega &= \sqrt{2} \left( -b^{157} + b^{236}\right) - b^4 \wedge \left(b^{16}-b^{27}-b^{35}\right) \\
g &= 2 \cdot (b^1 \cdot b^6 + b^2 \cdot b^7+ + b^3 \cdot b^5) - b^4 \cdot b^4
\end{align*}
and $\mathfrak{h}$ is a subalgebra of 
\begin{align*}
\mathfrak{h}^{I\!I} &{\coloneqq} \left\lbrace h(A,z,c) \, \middle|\, A \in \mathfrak{gl}(2,\mathbb{R}), \; z \in \mathbb{R}^{4}, c \in \mathbb{R} \right\rbrace = \mathfrak{p}_2 \hspace{.5em},
\end{align*}
where
\begin{align*}
h(A,z,c) = \begin{pmatrix}
A & \sigma(z) & U(c) \\ 0 & \rho(A) & \sigma(z)^* \\ 0 & 0 & -A^{\top} \end{pmatrix} \hspace{.5em}.
\end{align*}
\end{lemma}
See \cite{Fino2016} for a proof and most of the technical details about Type II holonomy algebras. Here, we concentrate on parts of the notation which are essential in order to understand the proof of Theorem \ref{Thm:MainTheorem}. First, we embed $\mathfrak{gl}(2,\mathbb{R})$ into $\mathfrak{g}_{2}^{*}$ by 
\begin{align*}
\mathfrak{gl}(2,\mathbb{R}) \cong \left\lbrace h(A,0,0) \,\middle|\, A \in \mathfrak{gl}(2,\mathbb{R}) \right\rbrace
\end{align*}
 and denote by $\mathfrak{a}$ the projection of $\mathfrak{h}\subset \mathfrak{h}^{I\!I}$ to $\mathfrak{gl}(2,\mathbb{R})$. Using the definitions
\begin{align*}
\mathfrak{n} &\coloneqq \left\lbrace h(0,z,c) \, \middle| \, z \in \mathbb{R}^4,\, c \in \mathbb{R} \right\rbrace \\
\mathfrak{n}(i,j) &\coloneqq \left\lbrace h(0,z,c) \, \middle| \, z \in \mathbb{R}^4,\, z_l = 0 \;\mathrm{if}\; l \notin \{i,j\},\, c \in \mathbb{R}
 \right\rbrace \\
 \mathfrak{n}(i,j,k) &\coloneqq \left\lbrace h(0,z,c) \, \middle| \, z \in \mathbb{R}^4,\, z_l = 0 \;\mathrm{if}\; l \notin \{i,j,k \},\, c \in \mathbb{R} \right\rbrace \hspace{.5em},
\end{align*}
where $i,j,k \in \{1,\ldots,4\}$, $\mathfrak{h}^{I\!I}$ is expressible as a semidirect product
\begin{align*}
\mathfrak{h}^{I\!I} = \mathfrak{gl}(2,\mathbb{R}) \ltimes \mathfrak{n} \hspace{.5em}.
\end{align*}
The action of $A \in \mathfrak{gl}(2,\mathbb{R})$ on $\mathfrak{n}$ is given by 
\begin{align*}
A \cdot h(0,z,c) = h(0,A \cdot z, \mathrm{tr}(A) \cdot c) \hspace{.5em}
\end{align*}
with respect to the representation of $\mathfrak{gl}(2,\mathbb{R})$ on $\mathbb{R}^{4}$ defined by the equation 
\begin{align*}
\sigma(A\cdot z ) = A \circ \sigma(z) - \sigma(z) \circ \rho(A) \hspace{.5em},
\end{align*}
where $\rho(A)$ and $\sigma(z)$ are the matrices defined above.

While the proof of Lemma \ref{Lem:TypeIIAlgebras} relies mainly on indecomposability and consequently on the dimension of the socle, we now want to determine which subalgebras of $\mathfrak{h}^{I\!I}$ are Berger algebras using Berger's first criterion. As we are not distinguish between locally symmetric manifolds and those who are not, we do not concern ourselves with Berger's second criterion.

Let $\mathcal{K}(\mathfrak{h})$ be the space of algebraic curvature endomorphisms with values in $\mathfrak{h}$ defined by 
\begin{align*}
\mathcal{K}(\mathfrak{h}) \coloneqq \left\lbrace R\in \Lambda^{2} V^{*} \otimes \mathfrak{h} \;\middle|\; \forall X,Y,Z \in V : \underset{X,Y,Z}{\mathfrak{S}} R(X,Y)Z = 0 \right\rbrace \hspace{.5em},
\end{align*}
where $\mathfrak{S}$ denotes summation over cyclic permutations . Now, consider the space 
\begin{align*}
\underline{\mathfrak{h}} \coloneqq \left\lbrace R(X,Y) \;\middle|\; X,Y \in V,\, R \in \mathcal{K}(\mathfrak{h}) \right\rbrace \subset \mathfrak{h} \hspace{.5em},
\end{align*}
and by Berger's criterion we have $\mathfrak{h} = \underline{\mathfrak{h}}$ for $\mathfrak{h}$ being a holonomy algebra. This implies for $R \in \mathcal{K}(\mathfrak{h})$ 
\begin{align*}
R_{ij}\coloneqq R(b_i,b_j) = h(A^{ij},z^{ij},c^{ij}) \hspace{.5em},
\end{align*}
with respect to the basis $b_1,\ldots,b_7$ chosen in Lemma \ref{Lem:TypeIIAlgebras}. Fino and Kath provide a parametrization of $\mathcal{K}(\mathfrak{h})$  \cite{Fino2016}, which is one of the key ingredients for the proof of both the classification theorem and Theorem \ref{Thm:MainTheorem}.
\begin{lemma}[\cite{Fino2016}]
The space $\mathcal{K}(\mathfrak{h})$ can be parametrized by real numbers $r,r_1,\ldots,r_4$, $j_1,j_2,t,t_1,\ldots,t_6,s_1,s_2,x_1,\ldots,x_5,y_1,\ldots,y_5$,  where $R = h(A,z,c) \in \mathcal{K}(\mathfrak{h})$ is given by the data in Table \ref{Tab:RII}.
\end{lemma}
\begin{table}
\begin{center}
\caption{Parametrization of curvature endomorphisms which are able to span holonomy algebras $\mathfrak{h}\subset \mathfrak{g}_2^*$. Adopted from \cite{Fino2016}. In the original table one has $r = 0$ which is wrong but has no considerable implication.}\label{Tab:RII}
\input{TabelleRTypII}

\end{center}
\end{table}
This parametrization together with some algebraic considerations about the structure of $\mathfrak{h} \subset \mathfrak{p}_2$ is consequently used in \cite{Fino2016} in order to proof the classification theorem for Type II algebras. Here, introducing an additional piece of notation is necessary before stating this theorem. Thus, denote by $C_a$ the matrix 
\begin{align*}
C_a \coloneqq \begin{pmatrix}
a & -1 \\ 1 & a
\end{pmatrix} \hspace{.5em},
\end{align*}
define algebras
\begin{align*}
Z_0&\coloneqq \left\lbrace h(0,(3z_1,z_2,z_1,3z_2),c) \,\middle|\, z_1,z_2,c \in \mathbb{R}\right\rbrace\\
Z_1&\coloneqq \left\lbrace h(0,(z_1,0,z_1,z_4),c)\;\middle|\; z_1,z_4,c \in \mathbb{R} \right\rbrace\\
Z_2&\coloneqq \left\lbrace h(0,(0,z_2,z_3,-z_2),c)\;\middle|\; z_2,z_3,c \in \mathbb{R} \right\rbrace
\end{align*}
and for parameters $\alpha, \beta, \kappa, s$ algebras
\begin{align*}
Z_3 = Z_3(\alpha)&\coloneqq \left\lbrace h(0,(z_1,\alpha z_1,\alpha z_4,z_4),c)\;\middle|\; z_1,z_4,c \in \mathbb{R} \right\rbrace \\
Z_4 = Z_4(s,\beta)&\coloneqq \left\lbrace h(0,(s z_1,-\beta z_4,-\beta z_1,z_4),c)\;\middle|\; z_1,z_4,c \in \mathbb{R} \right\rbrace \\
Z_5 = Z_5(\kappa)&\coloneqq \left\lbrace h(0,(z_1,z_2,\kappa z_1,z_4),c)\;\middle|\; z_1,z_2,z_4,c \in \mathbb{R} \right\rbrace \hspace{.5em},
\end{align*}
where
\begin{align*}
\kappa=\pm 1, \quad \alpha \in \left[ \frac{\sqrt{3}-1}{\sqrt{6}},\frac{\sqrt{3}+1}{\sqrt{6}}\right],\quad
s \in (0,1],\, \beta \in \mathbb{R}:\; 3\beta^2-(s+1)\beta -s =0 \hspace{.5em}.
\end{align*}
Furthermore, define
\begin{align*}
\mathfrak{d} \coloneqq \left\lbrace \mathrm{diag}(a,d) \, \middle|\, a,d \in \mathbb{R} \right\rbrace\qquad
\mathfrak{co}(2) \coloneqq \left\lbrace \begin{pmatrix} a & -b \\ b & a \end{pmatrix}\,\middle|\, a,b\in \mathbb{R}\right\rbrace \hspace{.5em}.
\end{align*}
Then, we have the following classification of possible holonomy algebras $\mathfrak{h}$ contained in $\mathfrak{gl}(2,\mathbb{R}) \ltimes \mathfrak{n} = \mathfrak{h}^{I\!I} \subset \mathfrak{g}_2^{*}$.
\begin{theorem}[Type II holonomy algebras, \cite{Fino2016}]\label{Thm:Classification}
If $\mathfrak{h}$ is of Type II, then there exists a basis of V such that we are in one of the following cases
\begin{itemize}
\item[(1)] $\mathfrak{a} \in \left\lbrace \mathfrak{gl}(2,\mathbb{R}),\mathfrak{sl}(2,\mathbb{R})\right\rbrace$ and $\mathfrak{h} = \mathfrak{a} \ltimes \mathfrak{n}$
\item[(2)] $\mathfrak{a}\in \left\lbrace \mathfrak{co}(2), \mathbb{R}\cdot C_a \right\rbrace$ and $\mathfrak{h} = \mathfrak{a} \ltimes \mathfrak{n}$ or $\mathfrak{h} = \mathfrak{a} \ltimes Z_0$.

\item[(3)] $\mathfrak{a} = \mathfrak{d}$ and $\mathfrak{h}= \mathfrak{d} \ltimes \mathfrak{n}_1$, where \vspace{.25ex}

{\centering
$\displaystyle
\mathfrak{n}_1 \in \left\lbrace \mathfrak{n},\mathfrak{n}(1,3),\mathfrak{n}(2,3),\mathfrak{n}(1,2,3),\mathfrak{n}(1,2,4) \right\rbrace$
\par}\vspace{.25ex}
\item[(4)] $\mathfrak{a} = \mathbb{R}\cdot \mathrm{diag}(1,\mu), \mu \in [-1,1)$, and 
\begin{itemize}
\item[(a)] $\mu \in [-1,1)$ and $\mathfrak{h} = \mathfrak{a} \ltimes \mathfrak{n}_1$, where \vspace{.25ex}

{\centering
$\displaystyle
\mathfrak{n}_1 \in \left\lbrace \mathfrak{n},\mathfrak{n}(2,3),\mathfrak{n}(1,2,3),\mathfrak{n}(1,2,4),\mathfrak{n}(1,3,4),\mathfrak{n}(2,3,4)\right\rbrace$
\par}\vspace{.25ex}
\item[(b)] $\mu=\frac{1}{2}$ and $\mathfrak{h}=\mathbb{R} \cdot h(\mathrm{diag}(1,\frac{1}{2}),(1,0,0,0),0) \ltimes \mathfrak{n}_1$, where \vspace{.25ex}

{\centering
 $\displaystyle
\mathfrak{n}_1 \in \left\lbrace \mathfrak{n}(2,3),\mathfrak{n}(2,3,4)\right\rbrace$
\par}\vspace{.25ex}
\item[(c)] $\mu = 0$ and $\mathfrak{h}= \mathfrak{a} \ltimes \mathfrak{n}(2,4)$ or \vspace{.25ex}

{\centering
 $\displaystyle
\mathfrak{h} = \mathbb{R}\cdot h(\mathrm{diag}(1,0),(0,1,0,0),0) \ltimes \mathfrak{n}_1 \hspace{.5em},$
\par}\vspace{.25ex} where $\mathfrak{n}_1 \in \left\lbrace \mathfrak{n}(1,4), \mathfrak{n}(3,4),\mathfrak{n}(1,3,4)\right\rbrace$
\end{itemize}
\item[(5)] $\mathfrak{a}\in \left\lbrace 0, \mathbb{R}\cdot I \right\rbrace$ and $\mathfrak{h} = \mathfrak{a} \ltimes \mathfrak{n}_1$, where \vspace{.25ex}

{\centering
$\displaystyle
\mathfrak{n}_1 \in \left\lbrace \mathfrak{n},\mathfrak{n}(1,3),\mathfrak{n}(2,3),\mathfrak{n}(1,2,4),\mathfrak{n}(2,3,4),Z_1,Z_2,Z_3,Z_4,Z_5 \right\rbrace \hspace{.5em}. $
\par}\vspace{.25ex}
\end{itemize}
\end{theorem} 
With this classification in mind, we review the procedure developed in \cite{Fino2017} to determine which of these Berger algebras are realizable by an appropriate metric in the next section.

\section{Local Type II metrics}
In this section we provide a construction scheme for metrics with certain holonomy properties. With this scheme it is possible to construct in each case of Theorem \ref{Thm:Classification} a metric such that the holonomy algebra of the associated Levi-Civita connection is equal to the particular Berger algebra. 

As before, let $M$ be a 1-connected, 7-dimensional manifold. Denote by $GL(M)$ the frame bundle over $M$. By the considerations in section \ref{Sec:2} we have an inclusion $\mathrm{G}_2^{*} \hookrightarrow SO(4,3) \hookrightarrow GL(7)$ and hence a reduction of the frame bundle to a principal $\mathrm{G}_2^{*}$-subbundle $\mathrm{G}_2^{*}(M)$. Via the 3-form $\omega$ a metric, the Levi-Civita connection $\nabla$ and consequently the associated connection form $\theta$  on $\mathrm{G}_2^{*}(M)$ are induced.

Now, we choose one of the Berger algebras $\mathfrak{h}$ from Theorem \ref{Thm:Classification} with Lie group $H$. Again, we have a reduction of  $\mathrm{G}_2^{*}(M)$ to a $H$-subbundle $H(M)$ along the inclusion $H \hookrightarrow \mathrm{G}_2^*$. The question now is if the Levi-Civita connection reduces to $H(M)$. It is well known from general holonomy theory that such a reduction always exists if $H$ is a holonomy group. A necessary and sufficient condition for the existence of such a reduction is that the connection form $\theta$ restricted to $H(M)$ takes values in $\mathfrak{h}$ only. 

We use this condition and the Ambrose-Singer holonomy theorem in order to construct in three subsequent steps a metric such that the Levi-Civita connection reduces to $H(M)$. Preliminary, for a given Berger algebra $\mathfrak{h}$, choose a local frame $(e_i) \in H(M)$ with local coframe $b = (b^i)$, $b: TH \to \mathbb{R}^7$. 

In the first step, \textit{assume} $H$ is a holonomy group. Under this assumption we deduce some necessary properties of the local coframe $b$. We start with Cartan's structure equation
\begin{align*}
\mathrm{d}b = - \theta \wedge b \hspace{.5em},
\end{align*}
where $\theta$ is restricted, by the holonomy condition, to take values in $\mathfrak{h}$ only. In components we have
\begin{align*}
\mathrm{d}b^i = -\sum \theta^{i}_{j} \wedge b^j \hspace{.5em},
\end{align*}
which gives rise to an exterior differential system. Thus, Cartan's methods dealing with such systems are applicable in order to solve for $b^i$. The components $\theta^{i}_{j}$ of $\theta$ are determined by the particular algebra $\mathfrak{h}$ chosen from Theorem \ref{Thm:Classification}.   Next, we choose local coordinates $x_1,\ldots,x_7$ on $M$ and use repeated applications of
\begin{enumerate}
\item transformations of the local frame by elements of $H$,
\item Frobenius' theorem
\end{enumerate}
in order to solve the system. We obtain expressions of the form
\begin{align*}
b^i = \mathrm{e}^{{f_i^i}} \mathrm{d}x_i + \sum_{j \neq i, j \geq 4} f_{i}^{j}\, \mathrm{d}x_j \hspace{.5em},
\end{align*}
where $f_{i}^{j}=f_{i}^{j}(x_{\alpha_1},\ldots,x_{\alpha_l})$, $\alpha_l,l \in \{1,\ldots,7\}$, are local functions on $M$. We call this new local coframe \textit{adapted} by the following reason: the procedure ensures that each local function depends on as few local coordinates as possible, i.e., $l$  is separately minimal for each function. Nevertheless, in general the above procedure does not yield an unique set of local coordinates such that the choice of an adapted coframe is also not unique.

Expressing functions via an exponential for $i=j$ has purely practical reasons as it becomes clear later on.

 Components $\theta^{i}_{j}$ are calculated using the adapted local coframe by
\begin{align*}
\theta^{i}_{j} = b^{i}(\nabla b_j)
\end{align*}
and Koszul's formula. Components have to fulfill certain relations among each other, arising from the structure of the algebra $\mathfrak{h}$. For instance, if we assume $\mathfrak{h} = \mathfrak{sl}(2,\mathbb{R}) \ltimes \mathfrak{n}$ we have $\theta^{1}_{1} \overset{!}{=} - \theta^{2}_{2}$ by Lemma \ref{Lem:TypeIIAlgebras}. These relations reformulate as a partial differential system for local functions $f_{i}^{j}$. Thus, our assumption that $H$ is a holonomy group with respect to the Levi-Civita connection, or equivalently that the Levi-Civita connection reduces to $H(M)$, is valid if there exists a solution of this partial differential system. Nevertheless, this is just a necessary condition. Using Ambrose-Singer theorem \cite{Ambrose1953} we finally need to show that the local functions $f_{i}^{j}$ can be chosen in a way such that the curvature endomorphisms
\begin{align*}
R_{ij} = R(b_i,b_j) \coloneqq \nabla_{b_i} \nabla _{b_j} - \nabla_{b_j} \nabla _{b_i} - \nabla_{[b_i,b_j]}
\end{align*}
 span $\mathfrak{h}$ with respect to the adapted coframe. In some cases it is necessary to consider covariant derivatives $\nabla R$  as well (cf. \cite{Ozeki1956}). 
 
In summary, the procedure consists of three steps.
 \begin{enumerate}
 \item Choose an algebra $\mathfrak{h}$ from Theorem \ref{Thm:Classification}. Obtain an exterior differential system by the structure equation. Solve for the $(b^i)$ and express them via local functions $f_{i}^{j}$.
 \item The structure of $\mathfrak{h}$ induces relations between the components of $\theta$. Use covariant differentiation in order to obtain a PDE system for the local functions $f_{i}^{j}$. Solve it as far as possible.
 \item Determine the $f_{i}^{j}$ such that they fulfill remaining equations from step 2 and $\mathfrak{h}$ is spanned by the curvature endomorphisms $R_{ij}$.
 \end{enumerate}
 Finally, if all of these calculations are performed successfully we obtain a metric $g$ by Eq. \eqref{Eq:InnerProduct} in terms of local functions on $M$ such that the manifold $(M,g)$ has holonomy $H \subset \mathrm{G}_2^{*}$ with respect to the induced Levi-Civita connection. 
 
The following example illustrates the described procedure. 

\begin{example}[Type II 5b.10: $\mathfrak{h} = \mathbb{R} \cdot I \ltimes Z_5$]
By the structure equation
\begin{align*}
\mathrm{d}b = - \theta \wedge b
\end{align*}
we deduce the exterior differential system
\begin{align*}
\mathrm{d}b^1 &= -\mathbf{x} \wedge b^1 -\mathbf{z}_2 \wedge b^3 - \sqrt{2}\kappa \mathbf{z}_1 \wedge b^4 - \mathbf{z}_4 \wedge b^5 + \mathbf{c} \wedge b^7\\
\mathrm{d}b^2 &= -\mathbf{x} \wedge b^2 -\mathbf{z}_1 \wedge b^3 - \sqrt{2} \mathbf{z}_2 \wedge b^4 -\kappa \mathbf{z}_1 \wedge b^5 - \mathbf{c} \wedge b^6\\
\mathrm{d}b^3 &=  \mathbf{z}_4 \wedge b^6  +\kappa \mathbf{z}_1 \wedge b^7\\
\mathrm{d}b^4 &= -\sqrt{2}\kappa \mathbf{z}_1 \wedge b^6  -\sqrt{2}\mathbf{z}_2 \wedge b^7\\
\mathrm{d}b^5 &= \mathbf{z}_2 \wedge b^6 + \mathbf{z}_1 \wedge b^7\\
\mathrm{d}b^6 &= \mathbf{x} \wedge b^6\\
\mathrm{d}b^7 &= \mathbf{x} \wedge b^7\hspace{.5em},
\end{align*}
where bold letters $\mathbf{x}$, $\mathbf{z}_i$, $\mathbf{c}$ denote real-valued differential 1-forms and $\kappa=\pm 1$. Using Frobenius' theorem and appropriate transformations of the local coframe by the action of the structure group $H$, we can solve this system in terms of local coordinates $x_1,\ldots,x_7$ on $M$. For a detailed calculation of these steps we refer to \cite{Volkhausen2019}. The techniques used there for Type III algebras are just as usable in case of Type II. We obtain
\begin{align*}
b^1 &= \mathrm{d}x_1 + r_6(x_3,x_4,x_5,x_6,x_7)\cdot \mathrm{d}x_6 + r_7(x_1,x_3,x_4,x_5,x_6,x_7)\cdot\mathrm{d}x_7\\
b^2 &= \mathrm{d}x_1 + s_7(x_2,x_3,x_4,x_5,x_6,x_7)\cdot \mathrm{d}x_7\\
b^3 &= \mathrm{d}x_3 + t_7(x_6,x_7)\cdot \mathrm{d}x_7\\
b^4 &= \mathrm{d}x_4 + u_6(x_6,x_7)\cdot \mathrm{d}x_6 + u_7(x_6,x_7)\cdot \mathrm{d}x_7\\
b^6 &= \mathrm{e}^{w_6(x_6,x_7)} \cdot \mathrm{d}x_6\\
b^5 &= \mathrm{d}x_5 \qquad b^7 = \mathrm{d}x_7 \hspace{.5em},
\end{align*} 
where we write $r_6,r_7,s_7,\ldots$ instead of $f_i^j$.
With the basis $b^i$ expressed in terms of local coordinates we have finished the first step. Next, the structure of $\mathfrak{h} = \mathbb{R} \cdot I \ltimes Z_5$ implies the following relations between components of the connection form $\theta$
\begin{align*}
\theta_{1}^{1} &= \theta_{2}^{2} \qquad \theta^{2}_{4} = \sqrt{2}\, \theta^{1}_{3} \qquad \theta^{1}_{4} = \sqrt{2}\, \theta^{2}_{5}\qquad \theta^{2}_{5} = \kappa \theta^{2}_{3} \hspace{.5em}.
\end{align*}
All other relations turn out to be either fulfilled trivially, by symmetry or contain no further information, i.e., they are redundant. Translating these relations into a PDE system, using $\theta^{i}_{j} = b^i(\nabla b_j)$ and Koszul's formula, we obtain
\begin{align*}
(r_7)_{x_1} &= (w_6)_{x_7} \\
(r_7)_{x_3} &=  \sqrt{2} (s_7)_{x_4} \\
(r_7)_{x_4} &= 2 \sqrt{2} (s_7)_{x_5} + \left[ (u_6)_{x_7} - (u_7)_{x_6} \right] \cdot \mathrm{e}^{-w_6} \\
(r_7)_{x_4} &= \left[ 2 \sqrt{2} (r_6)_{x_3} -  (u_6)_{x_7}  + (u_7)_{x_6} \right] \cdot \mathrm{e}^{-w_6}\\
(r_7)_{x_5} &= \kappa (r_7)_{x_3} - (t_7)_{x_6} \cdot \mathrm{e}^{-w_6}\\
(r_7)_{x_5} &= \left[ \sqrt{2} (r_6)_{x_4} - (t_7)_{x_6} \right] \cdot \mathrm{e}^{-w_6}\\
(s_7)_{x_2} &= (w_6)_{x_7}\\
(s_7)_{x_5} &= \kappa (s_7)_{x_3} \hspace{.5em}.
\end{align*}
Using $\mathrm{e}^{w_6}$ instead of $w_6$ avoid derivatives of the form $f^{-1}\cdot (f)_{x_i}$, which is quite useful in the more complicate cases. This is the `practical reason' mentioned above. The last equation allows for $s_7 \sim x_5 + \kappa^{-1} x_3$. Furthermore, the only difficulty of the system is the equation $(r_7)_{x_3} =  \sqrt{2} (s_7)_{x_4}$.  With the ansatz $s_7 \sim f(x_6,x_7) \cdot \left( x_5 + \kappa^{-1} x_3 \right) + G_7'(x_2,x_6,x_7)$ the system is solvable by straightforward substitution and integration. We obtain
\begin{align*}
r_6 &= x_3 \cdot f(x_6,x_7) \cdot \mathrm{e}^{w_6}+ \frac{x_3}{\sqrt{2}} \left[ (u_6)_{x_7} - (u_7)_{x_6} \right]  + F_6(x_5,x_6,x_7)\\
r_7 &=x_1 (w_6)_{x_7} + 2\sqrt{2} x_4 f(x_6,x_7)+\left[ x_4 \left((u_6)_{x_7} - (u_7)_{x_6}\right) - x_5 (t_7)_{x_6} \right] \mathrm{e}^{-w_6} \\
&\hspace{9.5cm}+ F_7(x_6,x_7)\\
s_7 &= x_2 (w_6)_{x_7} + f(x_6,x_7) \cdot \left( x_5 + \kappa^{-1} x_3 \right) + G_7(x_6,x_7)\hspace{.5em},
\end{align*}
where $t_7,u_6,u_7$ are arbitrary functions up to this point. Note, the solution is not unique. Determining these functions such that the curvature endomorphisms span the algebra $\mathfrak{h}=\mathbb{R} \cdot I \ltimes Z_5$ is the task of the final step. With the choice
\begin{align*}
f &= \kappa x_6 - x_6 x_7 \qquad w_6 = \frac{1}{2} x_6 x_7^2 \qquad F_6 = \frac{x_5^2}{2}\\
&\hspace{1.5em}u_6 = u_7 = t_7 = F_7 = G_7 =0
\end{align*}
it is simple to calculate that the curvature endomorphisms $R_{ij}$ are linear combinations of 
\begin{align*}
R_{37} = h(0,0,1) \qquad R_{56} = h(0,(0,0,0,1),0) \qquad R_{67} = h(0,(1,0,\kappa,0),0) \\
\nabla_{b_6}R_{67} = h(0,(0,1,0,0),0) \qquad \nabla_{b_7}R_{67} = h\left(\mathrm{diag}(1,1),(-\frac{1}{\kappa},0,-1,0),0\right)
\end{align*}
and therefore span the algebra $\mathfrak{h}$. Consequently, the subgroup $H \subset \mathrm{G}_2^{*}$ is the holonomy group associated with the metric given by Eq. \eqref{Eq:InnerProduct}, the coframe $b^i$ and the chosen local functions.   
\end{example}

With the techniques used in the example we can prove Theorem \ref{Thm:MainTheorem} by analogue calculations in all cases of Theorem \ref{Thm:Classification}.
\begin{proof}[Main Theorem]
In each case of Theorem \ref{Thm:Classification} the basis and local functions listed in the appendix realize a metric $g$ such that the holonomy algebra with respect to the induced Levi-Civita connection is equal to the respective algebra in the theorem.
\end{proof}
\section*{Appendix}
The appendix is to read as follows: Table \ref{Tab:Normalforms} provides  how the basis $b^i$ is expressed in terms of local functions $f^j_{i}$ and  coordinates $x_i$ on $M$. Instead of using the general expression 
\begin{align*}
b^i = \mathrm{e}^{f_i^i}\mathrm{d}x_i + \sum_{j \neq i, j \geq 4} f^{j}_{i} \mathrm{d}x_j
\end{align*}
we assign a certain letter to local functions of each $b^i$ as in the example. Thus, we express $b^1,\ldots,b^7$ as 
\begin{align*}
b^1 &=  \mathrm{d}x_1+ \sum_{ i \geq 4} r_i \, \mathrm{d}x_i\;,\qquad
b^2 =  \mathrm{d}x_2 + \sum_{i \geq 4} s_i \, \mathrm{d}x_i\;,\qquad
b^3 =  \mathrm{d}x_3 +  \sum_{i \geq 6} t_i  \,\mathrm{d}x_i\;,\\
b^4 &=  \mathrm{d}x_4 + \sum_{i \geq 6} u_i \, \mathrm{d}x_i\;, \qquad
b^5 =  \mathrm{d}x_5 + \sum_{i \geq 6} v_i \, \mathrm{d}x_i\;,\qquad
b^6 = \mathrm{e}^{w_6}\,\mathrm{d}x_6 + w_7\, \mathrm{d}x_7\;,\\
b^7 &=   \mathrm{e}^{z_7}\,\mathrm{d}x_7 \hspace{.5em},
\end{align*}
where $i \in \{4,\ldots,7\}$ and $r_i,s_i,t_i,u_i,v_i,w_i,z_i$ are local functions. Note, there are only a few cases where $w_6, z_7 \neq 0$ and we dropped all functions which are not needed further.

In the following subsections local functions which generate the particular holonomy algebra  in each case of Theorem \ref{Thm:Classification} are listed. Since dependencies on local coordinates vary from case to case, they are not listed. Instead, they can be extracted in each case from Table \ref{Tab:Normalforms}.

Since the calculations themselves are quite lengthy and somewhat repetitive, they are not given here. Readers interested in details can obtain the Maple files used to compute the holonomy algebras via the author's website and on request. In case of any questions or problems with these files do not hesitate to contact the author.

\footnotesize
\begin{landscape}
\input{Normalforms_hor}
\end{landscape}
\normalsize
\input{Table_Type_II_1}
\input{Table_Type_II_2}
\input{Table_Type_II_3}
\input{Table_Type_II_4}
\input{Table_Type_II_5}
\clearpage
\bibliographystyle{amsplain}
\bibliography{Bibliography_Local_Type_II_Metrics}

\end{document}

%% file: TabelleRTypII.tex
\begin{tabular}{|c|c|c|c|}
\hline &&& \\[-2ex]
$R_{ij}=R\!\left(b_i,b_j \right)$ & $A$ & $z$ & $c$ \\[.75ex]
\hline &&&\\[-1.5ex]
$R_{16}$ & $0$ & $(x_4,x_3,x_2,x_1)$  & $t+t_1$ \\[1ex]
$R_{17}=\frac{1}{\sqrt{2}}R_{34}$ & $0$ & $(x_5,x_4,x_3,x_2)$  & $t_4-t_5$ \\[1ex]
$R_{26}=-\frac{1}{\sqrt{2}}R_{45}$ & 0 & $(y_4,y_3,y_2,y_1)$  & $t_2-t_3$ \\[1ex]
$R_{27}$ & $0$ & $(y_5,y_4,y_3,y_2)$  & $t_6+t$ \\[1ex]

$R_{56}$  &$\begin{pmatrix}x_1 & y_1 \\ x_2 & y_2\\ \end{pmatrix}$ & $(t_6,t_2,s_2,j_2)$  & $r_1$ \\[3ex]
$R_{57}=\frac{1}{\sqrt{2}}R_{46}$  &$\begin{pmatrix} x_2 & y_2 \\x_3 & y_3\\ \end{pmatrix}$ & $(t_5,t_1,t_3,s_2)$  & $r_2$ \\[2ex]
$R_{36}=\frac{1}{\sqrt{2}}R_{47}$ & $\begin{pmatrix} x_3 & y_3 \\ x_4 & y_4\\ \end{pmatrix}$ & $(s_1,t_4,t_1,t_2)$  & $r_3$ \\[3ex]
$R_{37} $ & $\begin{pmatrix} x_4 & y_4 \\ x_5 &y_5\\ \end{pmatrix}$ & $(j_1,s_1,t_5,t_6)$  & $r_4$ \\[3ex]
$R_{67}$& $\begin{pmatrix} t+t_1 & t_2-t_3 \\ t_4-t_5 &t_6+t\\ \end{pmatrix}$ & $(r_4,r_3,r_2,r_1)$  & $r$ \\[2ex]
\hline
\end{tabular}

%% file: Normalforms_hor.tex
\setlength{\LTcapwidth}{.8\textwidth}
\begin{longtable}[c]{llcl}
\caption{Local functions used for the adapted coframes $(b^i)$. Functions not listed are identically zero. }\label{Tab:Normalforms}\\
Nr. &\phantom{Alg} $\mathfrak{h}$ & $\dim \mathfrak{h}$ & Local functions \\\toprule \endfirsthead
\caption{Local functions used for (quasi-)normalforms of the basis $b^i$. Functions not listed are identically zero.} \\
Nr. &\phantom{Alg} $\mathfrak{h}$ & $\dim \mathfrak{h}$ & Local functions \\\toprule \endhead
\bottomrule\multicolumn{4}{r}{ \textit{Continued on next page} $\rightarrow$}\\\endfoot
\bottomrule\multicolumn{4}{r}{}\endlastfoot 
1a  & $\mathfrak{gl}(2,\mathbb{R})\ltimes \mathfrak{n}$ & 9 & $r_4(x_3,\ldots,x_7), r_5(x_3,\ldots,x_7), r_6(x_1,\ldots,x_7),r_7(x_1,\ldots,x_7), s_4(x_3,\ldots,x_7), s_5(x_3,\ldots,x_7),s_7(x_1,\ldots,x_7),$\\
&&& $t_7(x_3,\ldots,x_7), u_7(x_3,\ldots,x_7)$ \\
1b & $\mathfrak{sl}(2,\mathbb{R})\ltimes \mathfrak{n}$ & 8 & $r_4(x_3,\ldots,x_7), r_5(x_3,\ldots,x_7), r_6(x_1,\ldots,x_7),r_7(x_1,\ldots,x_7), s_4(x_3,\ldots,x_7), s_5(x_3,\ldots,x_7),s_7(x_1,\ldots,x_7),$\\
&&& $ t_7(x_3,\ldots,x_7), u_7(x_3,\ldots,x_7)$ \\\midrule
2a & $\mathfrak{co}(2)\ltimes \mathfrak{n}$ & 7 & $r_4(x_3,\ldots,x_7),r_5(x_3,\ldots,x_7), r_6(x_1,\ldots,x_7),r_7(x_1,\ldots,x_7), s_4(x_3,\ldots,x_7), s_5(x_3,\ldots,x_7), s_7(x_1,\ldots,x_7), t_7(x_3,\ldots,x_7),$\\
&&& $ u_7(x_3,\ldots,x_7), w_7(x_6,x_7)$\\
2b &$\mathfrak{co}(2)\ltimes Z_0$ & 5 & $r_4(x_3,\ldots,x_7), r_6(x_1,\ldots,x_7), r_7(x_1,\ldots,x_7), s_4(x_3,\ldots,x_7), s_7(x_1,\ldots,x_7),t_6(x_4,x_6,x_7),t_7(x_4,x_6,x_7),u_6(x_3,\ldots,x_7),$\\
&&& $ u_7(x_3,\ldots,x_7), w_7(x_6,x_7)$\\
2c &$\mathbb{R}\cdot C_a\ltimes \mathfrak{n}$ & 6 & $r_4(x_3,\ldots,x_7), r_6(x_1,\ldots,x_7),r_7(x_1,\ldots,x_7),s_4(x_3,x_4,x_6,x_7), s_5(x_3,x_5,x_6,x_7), s_7(x_1,\ldots,x_7), t_7(x_3,\ldots,x_7), $\\
&&& $u_7(x_3,\ldots,x_7), w_7(x_6,x_7), z_7(x_6,x_7)$\\
2d &$\mathbb{R}\cdot C_a\ltimes Z_0$ & 4 &  $r_4(x_3,\ldots,x_7), r_6(x_1,\ldots,x_7),r_7(x_1,\ldots,x_7),s_4(x_3,\ldots,x_7), s_7(x_1,\ldots,x_7),t_6(x_4,x_6,x_7),t_7(x_4,x_6,x_7),u_6(x_3,\ldots,x_7), $\\
&&& $u_7(x_3,\ldots,x_7), w_7(x_6,x_7), z_7(x_6,x_7)$\\\midrule
3a &  $\mathfrak{d}\ltimes \mathfrak{n}$& 7 & $r_5(x_3,x_5,x_6,x_7), r_6(x_1,x_3,\ldots,x_7),r_7(x_3,\ldots,x_7),s_5(x_3,\ldots,x_7), s_7(x_2,\ldots,x_7), t_7(x_3,x_6,x_7), u_7(x_5,x_6,x_7)$\\
3b & $\mathfrak{d}\ltimes \mathfrak{n}(1,3)$ & 5 & $r_6(x_1,x_4,x_6,x_7), r_7(x_3,x_5,x_6,x_7), s_5(x_3,x_5,x_6,x_7), s_7(x_2,x_3,x_5,x_6,x_7), t_6(x_3,x_6,x_7), t_7(x_3,x_6,x_7), v_6(x_5,x_6,x_7)$ \\
3c & $\mathfrak{d}\ltimes \mathfrak{n}(2,3)$& 5 & $r_5(x_3,x_5,x_6,x_7), r_6(x_1,x_3,\ldots,x_7), r_7(x_3,\ldots,x_7),s_5(x_4,\ldots,x_7), s_7(x_2,x_4,\ldots,x_7), t_6(x_3,x_6,x_7), t_7(x_3,x_6,x_7),$\\
&&& $ u_7(x_5,x_6,x_7), v_7(x_5,x_6,x_7)$ \\
3d & $\mathfrak{d}\ltimes \mathfrak{n}(1,2,3)$& 6 & $r_5(x_3,x_5,x_6,x_7), r_6(x_1,x_3,\ldots,x_7), r_7(x_3,\ldots,x_7),s_5(x_3,\ldots,x_7), s_7(x_2,\ldots,x_7), t_6(x_3,x_6,x_7), t_7(x_3,x_6,x_7), u_7(x_5,x_6,x_7)$\\
3e & $\mathfrak{d}\ltimes \mathfrak{n}(1,2,4)$& 6 & $r_5(x_3,x_5,x_6,x_7), r_6(x_1,x_3,x_5,x_6,x_7), r_7(x_3,\ldots,x_7), s_5(x_3,\ldots,x_7), s_7(x_2,\ldots,x_7), t_7(x_3,x_6,x_7), u_7(x_5,x_6,x_7),$\\\midrule
\multicolumn{4}{l}{ $\mathfrak{a}= \mathbb{R}\cdot\mathrm{diag}(1,\mu), \mu \in [-1,1)$; $Y_0 = \mathbb{R}\cdot h(\mathrm{diag(1,1/2),(1,0,0,0),0)}$} \\\midrule
4a.1 & $\mathfrak{a}\ltimes \mathfrak{n}$& 6 & $r_6(x_3,\ldots,x_7),r_7(x_1,x_3,\ldots,x_7),s_5(x_3,x_5,x_6,x_7), s_7(x_2,\ldots,x_7), t_7(x_3,x_6,x_7), u_7(x_6,x_7), v_7(x_5,x_6,x_7), w_6(x_6,x_7)$\\
4a.2 &$\mathfrak{a}\ltimes \mathfrak{n}(2,3)$& 4 & $r_6(x_3,x_4,x_6,x_7),r_7(x_1,x_3,\ldots,x_7),s_7(x_2,x_4,\ldots,x_7), t_7(x_3,x_6,x_7), u_7(x_6,x_7), v_7(x_5,x_6,x_7), w_6(x_6,x_7)$\\
4a.3 &$\mathfrak{a}\ltimes \mathfrak{n}(1,2,3)$& 5 & $r_6(x_3,x_4,x_6,x_7),r_7(x_1,x_3,\ldots,x_7),s_5(x_3,x_5,x_6,x_7), s_7(x_2,\ldots,x_7), t_7(x_3,x_6,x_7), u_7(x_6,x_7), w_6(x_6,x_7)$\\
4a.4 &$\mathfrak{a}\ltimes \mathfrak{n}(1,2,4)$& 5 & $r_6(x_3,x_5,x_6,x_7),r_7(x_1,x_3,x_5,x_6,x_7),s_5(x_3,x_5,x_6,x_7), s_7(x_2,\ldots,x_7), t_7(x_3,x_6,x_7), u_7(x_6,x_7), w_6(x_6,x_7)$\\ 
4a.5 &$\mathfrak{a}\ltimes \mathfrak{n}(1,3,4)$& 5 & $r_6(x_4,x_5,x_6,x_7),r_7(x_1,x_5,x_6,x_7),s_5(x_3,x_5,x_6,x_7), s_7(x_2,x_3,x_5,x_6,x_7), t_7(x_3,x_6,x_7), w_6(x_6,x_7)$\\
4a.6 &$\mathfrak{a}\ltimes \mathfrak{n}(2,3,4)$& 5 & $r_6(x_3,\ldots,x_7),r_7(x_1,x_3,\ldots,x_7),s_7(x_2,x_4,\ldots,x_7),  t_7(x_3,x_6,x_7), u_7(x_6,x_7), v_7(x_5,x_6,x_7), w_6(x_6,x_7)$\\
4b.1 &$Y_0\ltimes \mathfrak{n}(2,3)$ & 4 & $r_6(x_3,x_4,x_6,x_7),r_7(x_1,x_3,\ldots,x_7),s_7(x_2,\ldots,x_7), t_7(x_3,x_6,x_7), u_7(x_6,x_7), v_7(x_5,x_6,x_7), w_6(x_6,x_7)$\\
4b.2 &$Y_0\ltimes \mathfrak{n}(2,3,4)$& 5 & $r_6(x_3,\ldots,x_7),r_7(x_1,x_3,\ldots,x_7),s_7(x_2,\ldots,x_7), t_7(x_3,x_6,x_7), u_7(x_6,x_7), v_7(x_5,x_6,x_7), w_6(x_6,x_7)$\\\pagebreak
\multicolumn{4}{l}{ $\mathfrak{a}= \mathbb{R}\cdot\mathrm{diag}(1,0)$; $Y_1 = \mathbb{R}\cdot h(\mathrm{diag(1,0),(0,1,0,0),0)}$ } \\\midrule
4c.1 &$\mathfrak{a}\ltimes \mathfrak{n}(2,4)$& 4 & $r_5(x_3,x_5,x_6,x_7), r_6(x_1,x_3,x_5,x_6,x_7),r_7(x_4,\ldots,x_7), s_5(x_4,\ldots,x_7), u_7(x_5,x_6,x_7)$\\
4c.2 &$Y_1\ltimes \mathfrak{n}(1,4)$& 4 & $r_5(x_3,x_5,x_6,x_7), r_6(x_1,x_3,x_5,x_6,x_7),r_7(x_3,\ldots,x_7), s_7(x_3,x_6,x_7), u_7(x_6,x_7), v_6(x_5,x_6,x_7)$\\
4c.3 &$Y_1\ltimes \mathfrak{n}(3,4)$& 4 & $r_5(x_3,\ldots,x_7), r_6(x_1,x_3,\ldots,x_7),r_7(x_4,\ldots,x_7), s_7(x_5,x_6,x_7), t_7(x_6,x_7), v_6(x_5,x_6)$\\
4c.4 &$Y_1\ltimes \mathfrak{n}(1,3,4)$& 5 & $r_5(x_3,x_5,x_6,x_7), r_6(x_1,x_3,\ldots,x_7),r_7(x_3,\ldots,x_7), s_7(x_3,\ldots,x_7), u_6(x_5,x_6,x_7), u_7(x_5,x_6,x_7), v_6(x_5,x_6,x_7)$\\\midrule
5a.1 &  $ \mathfrak{n}$ & 5 & $r_6(x_3,\ldots,x_7), r_7(x_3,\ldots,x_7),s_7(x_3,\ldots,x_7),u_6(x_6,x_7), u_7(x_6,x_7)$\\
5a.2 & $\mathfrak{n}(1,3)$ & 3 & $r_6(x_4,x_6,x_7),r_7(x_5,x_6,x_7),s_7(x_3,x_6,x_7),t_7(x_6,x_7)$ \\
5a.3 & $\mathfrak{n}(2,3)$ & 3 & $r_6(x_3,x_6,x_7),r_7(x_3,x_4,x_6,x_7),s_7(x_5,x_6,x_7),u_7(x_6,x_7)$\\
5a.4 & $\mathfrak{n}(1,3,4)$ & 4 & $r_6(x_4,\ldots,x_7),r_7(x_5,x_6,x_7),s_7(x_3,x_6,x_7),t_7(x_6,x_7)$\\
5a.5 & $\mathfrak{n}(2,3,4)$ & 4 & $r_6(x_3,\ldots,x_7),r_7(x_4,\ldots,x_7),s_5(x_5,x_6,x_7), u_6(x_6,x_7),u_7(x_6,x_7)$\\
5a.6 & $Z_1$ & 3 & $ r_6(x_5,x_6,x_7),r_7(x_4,\ldots,x_7),s_7(x_3,x_5,x_6,x_7),t_7(x_6,x_7), u_7(x_6,x_7)$\\
5a.7 & $Z_2$ & 3 & $r_6(x_3,\ldots,x_7),r_7(x_4,\ldots,x_7),s_7(x_5,x_6,x_7),t_6(x_6,x_7),t_7(x_6,x_7)$\\
5a.8 & $Z_3$ & 3 & $r_6(x_4,\ldots,x_7),r_7(x_3,\ldots,x_7),s_7(x_3,\ldots,x_7),t_6(x_6,x_7), t_7(x_6,x_7), u_6(x_6,x_7), v_7(x_6,x_7)$\\
5a.9 & $Z_4$ & 3 &  $r_6(x_4,\ldots,x_7),r_7(x_3,\ldots,x_7),s_7(x_3,\ldots,x_7),t_6(x_6,x_7), t_7(x_6,x_7), u_6(x_6,x_7), u_7(x_6,x_7)$\\
5a.10 & $Z_5$ & 4 & $r_6(x_4,\ldots,x_7),r_7(x_3,\ldots,x_7),s_7(x_3,\ldots,x_7),t_7(x_6,x_7), u_6(x_6,x_7), u_7(x_6,x_7)$\\
5b.1 & $\mathbb{R}\cdot I \ltimes\mathfrak{n}$ & 6 & $r_6(x_1,x_3,\ldots,x_7),r_7(x_1,x_3,\ldots,x_7),s_7(x_2,\ldots,x_7), u_6(x_6,x_7), u_7(x_6,x_7), w_6(x_6,x_7)$\\
5b.2 & $\mathbb{R}\cdot I \ltimes\mathfrak{n}(1,3)$ & 4 & $s_6(x_1,x_4,x_6,x_7),s_7(x_2,\ldots,x_7),u_7(x_6,x_7),w_6(x_6,x_7)$\\
5b.3 & $\mathbb{R}\cdot I \ltimes\mathfrak{n}(2,3)$ & 4 & $r_6(x_3,x_4,x_6,x_7),r_7(x_1,x_4,x_5,x_6,x_7),s_7(x_2,x_5,x_6,x_7), u_7(x_6,x_7), w_6(x_6,x_7)$\\
5b.4 & $\mathbb{R}\cdot I \ltimes\mathfrak{n}(1,3,4)$ & 5 & $r_6(x_4,\ldots,x_7),r_7(x_1,x_5,x_6,x_7),s_7(x_2,x_3,x_5,x_6,x_7), t_7(x_6,x_7), w_6(x_6,x_7)$\\
5b.5 & $\mathbb{R}\cdot I \ltimes\mathfrak{n}(2,3,4)$ & 5 & $r_6(x_3,\ldots,x_7),r_7(x_1,x_4,\ldots,x_7),s_7(x_2,x_5,x_6,x_7),u_6(x_6,x_7), w_6(x_6,x_7)$\\
5b.6 & $\mathbb{R}\cdot I \ltimes Z_1$ & 4 &$ r_6(x_5,x_6,x_7),r_7(x_1,x_4 ,\ldots,x_7),s_7(x_2,x_3,x_5,x_6,x_7), t_7(x_6,x_7),u_7(x_6,x_7),w_6(x_6,x_7)$\\
5b.7 & $\mathbb{R}\cdot I \ltimes Z_2$ & 4 & $r_6(x_3,\ldots,x_7),r_7(x_1,x_4,\ldots,x_7), s_7(x_2,x_5,x_6,x_7),  t_6(x_6,x_7), t_7(x_6,x_7),w_6(x_6,x_7)$\\
5b.8 & $\mathbb{R}\cdot I \ltimes Z_3$ & 4 & $r_6(x_3,\ldots,x_7),r_7(x_1,x_3,\ldots,x_7),s_7(x_2,\ldots,x_7), t_6(x_6,x_7), t_7(x_6,x_7), u_6(x_6,x_7), v_7(x_6,x_7), w_6(x_6,x_7)$\\
5b.9 & $\mathbb{R}\cdot I \ltimes Z_4$ & 4 & $r_6(x_3,\ldots,x_7),r_7(x_1,x_3,\ldots,x_7),s_7(x_2,\ldots,x_7), t_6(x_6,x_7), t_7(x_6,x_7), u_6(x_6,x_7), u_7(x_6,x_7), w_6(x_6,x_7)$\\
5b.10 & $\mathbb{R}\cdot I \ltimes Z_5$ & 5 & $r_6(x_3,\ldots,x_7),r_7(x_1,x_3,\ldots,x_7),s_7(x_2,\ldots,x_7),t_7(x_6,x_7), u_6(x_6,x_7),u_7(x_6,x_7),w_6(x_6,x_7)$\\
\end{longtable}

%% file: Table_Type_II_1.tex
\clearpage
\begin{center}
\textsc{Type II 1}
\end{center}
Let $\mathfrak{a} \in \{\mathfrak{gl}(2,\mathbb{R}), \mathfrak{sl}(2,\mathbb{R})\}$ and $\mathfrak{h} = \mathfrak{a}\ltimes \mathfrak{n}$.
\small
\settowidth{\temp}{$\quad+2\sqrt{2}x_3  x_5 \left(x_4\left( \sqrt{2} x_4  -  x_7 \sqrt{2} \right)-\sqrt{2}x_3 x_5  \right)$}
\begin{align*}
r_4 &= \begin{cases}
\makebox[\temp][l]{$2\sqrt {2} x_3 x_5+4\sqrt {2} x_4^{2}-4 x_4 x_7 \,,$} & \mathrm{if\;} \mathfrak{a} = \mathfrak{gl}(2,\mathbb{R})\\[.25ex]
\sqrt{2}\left( x_5-\frac{1}{2} x_5^{2} \right) \,, & \mathrm{if\;} \mathfrak{a} = \mathfrak{sl}(2,\mathbb{R})\\[.25ex]
\end{cases}\\[1ex]
r_5 &= \begin{cases}
\makebox[\temp][l]{$2 x_3 ( x_7	-\sqrt{2} x_4) \,,$} & \mathrm{if\;} \mathfrak{a} = \mathfrak{gl}(2,\mathbb{R})\\[.25ex]
-3\sqrt{2} x_4 \left(x_5-1 \right) \,,	& \mathrm{if\;} \mathfrak{a} = \mathfrak{sl}(2,\mathbb{R})\\[.25ex]
\end{cases}\\[1ex]
r_6 &= \begin{cases}
\makebox[\temp][l]{$2 x_2 x_3 +  \sqrt{2}x_1 x_4 + 2  x_4^{2} x_7^{2}$}\\\quad
 +2\sqrt{2}x_3 x_5\left(  \sqrt{2} \left(  x_4^{2}-   x_3 x_5\right) - x_4  x_7 \right)\\\qquad
 +\frac{16}{3} x_4^{3} \left( x_4- \sqrt{2} x_7 \right)  \,,
	& \mathrm{if\;} \mathfrak{a} = \mathfrak{gl}(2,\mathbb{R})\\[.25ex]
 x_2\left( x_5-1\right)+\sqrt{2}x_4\left[ x_4^{2}\left( x_3- 1\right)\right.\\\qquad
 \left.+\frac{3}{2} \left( x_3^{2}  -2 x_3   -x_5^{2}\right) \left( x_5-1\right)\right] \,,
	& \mathrm{if\;} \mathfrak{a} = \mathfrak{sl}(2,\mathbb{R})\\[.25ex]
\end{cases}\\[1ex]
r_7 &= \begin{cases}
\makebox[\temp][l]{$2 x_1 x_5+2 x_2\left( 2 \sqrt{2} x_4-x_7\right) - \sqrt{2}  x_4x_7$}\\\quad
 + 4x_3 x_5^{2}\left(  x_7-2 \sqrt{2} x_4 \right)\\\qquad
 +2\sqrt{2}x_4 x_5\left(\frac{14}{3} x_4^{2} -4\sqrt{2} x_4 x_7 + x_7^{2} \right) \,,
	& \mathrm{if\;} \mathfrak{a} = \mathfrak{gl}(2,\mathbb{R})\\[.25ex]
	6x_4^{2}\left(   x_3 x_5 -   x_3-  x_5 +1 \right)+x_4^{4} \,,
	& \mathrm{if\;} \mathfrak{a} = \mathfrak{sl}(2,\mathbb{R})\\[.25ex]
\end{cases}\\[1ex]
s_4 &= \begin{cases}
\makebox[\temp][l]{$\frac{1}{\sqrt{2}}+2 x_4 x_5  -\sqrt{2} x_5  x_7 \,,$}
	& \mathrm{if\;} \mathfrak{a} = \mathfrak{gl}(2,\mathbb{R})\\[.25ex]
\frac{1}{\sqrt{2}}\left( x_3^{2}- 2x_3  -3  x_5^{2} \right) \,,
	& \mathrm{if\;} \mathfrak{a} = \mathfrak{sl}(2,\mathbb{R})\\[.25ex]
\end{cases}\\[1ex]
s_5 &= \begin{cases}
\makebox[\temp][l]{$0  \,, $}
	& \mathrm{if\;} \mathfrak{a} = \mathfrak{gl}(2,\mathbb{R})\\[.25ex]
3 x_5\left(1- \sqrt{2} x_4\right)	 \,,
	& \mathrm{if\;} \mathfrak{a} = \mathfrak{sl}(2,\mathbb{R})\\[.25ex]
\end{cases}\\[1ex]
s_7 &= \begin{cases}
\makebox[\temp][l]{$2 x_2 x_5  +\frac{1}{2} x_3^{2}-\frac{8}{3}  x_3 x_5^{3} \,,$} 	& \mathrm{if\;} \mathfrak{a} = \mathfrak{gl}(2,\mathbb{R})\\[.25ex]
 x_1 \left( x_3-1\right)  +\sqrt{2} x_4^{3} \left(x_5-1 \right) \,,
 	& \mathrm{if\;} \mathfrak{a} = \mathfrak{sl}(2,\mathbb{R})\\[.25ex]
\end{cases}\\[1ex]
t_7 &= \begin{cases}
\makebox[\temp][l]{$-2\left( x_3 x_5 + x_4^{2} \right)+\sqrt{2} x_4 x_7 \,,$}	& \mathrm{if\;} \mathfrak{a} = \mathfrak{gl}(2,\mathbb{R})\\[.25ex]
 x_5-\frac{1}{2} x_5^{2} \,,	& \mathrm{if\;} \mathfrak{a} = \mathfrak{sl}(2,\mathbb{R})\\[.25ex]
\end{cases}\\[1ex]
u_7 &= \begin{cases}
\makebox[\temp][l]{$-2 x_4 x_5 \,,$}	& \mathrm{if\;} \mathfrak{a} = \mathfrak{gl}(2,\mathbb{R})\\[.25ex]
\sqrt{2}  \left(x_3-\frac{1}{2} x_3^{2} \right) \,, & \mathrm{if\;} \mathfrak{a} = \mathfrak{sl}(2,\mathbb{R})\\[.25ex]
\end{cases}\\[1ex]
\end{align*}
\normalsize

%% file: Table_Type_II_2.tex
\begin{center}
\textsc{Type II 2}
\end{center}
Let $\mathfrak{a}\in \left\lbrace \mathfrak{co}(2),\mathbb{R}\cdot C_a \right\rbrace$ and $\mathfrak{h} = \mathfrak{a}\ltimes \mathfrak{n}$ or $\mathfrak{h}=\mathfrak{a}\ltimes Z_0$.
\small
%
%
\begin{align*}
r_4 &= \begin{cases}
\makebox[\temp][l]{$-  x_4 \left(\frac{1}{2} x_6^{2}+{ x_7} \right)  \,, $}
	& \mathrm{if\;} \mathfrak{h} = \mathfrak{co}(2) \ltimes \mathfrak{n}\\[.25ex]
 \left( x_7 - \frac{20}{3} \right) x_4 x_6 \,,
	& \mathrm{if\;} \mathfrak{h} = \mathfrak{co}(2) \ltimes Z_0\\[.25ex]
2\sqrt{2}  \left( a \left( x_3- x_5\right)- \sqrt{2} x_4 \right) x_6 \,,
	& \mathrm{if\;} \mathfrak{h} = \mathbb{R}\cdot C_a \ltimes \mathfrak{n}\\[.25ex]
-\frac {\sqrt{2} a \left(  x_3+3 x_5 \right) }{ a^{2}+1} \,,
	& \mathrm{if\;} \mathfrak{h} = \mathbb{R}\cdot C_a \ltimes Z_0\\[.25ex]
\end{cases}\\[1ex]
%
%
r_5 &= \makebox[\temp][l]{$ 0 $}	\qquad \mathrm{in\;all\;cases}\\
\end{align*}%
\settowidth{\temp}{$\quad -2 x_1 x_6  - x_2 x_6^{2} + \frac{9}{2} x_3^{2} + \left( 3 x_5 + \frac{8 \sqrt{2}}{3} x_4 x_6^{3} \right) x_3 + \frac{1}{2} x_5^{2} + \frac{8\sqrt{2}}{9} x_4 x_5 x_6^{3}  \,,$}%
%
%
\begin{align*}
r_6 &= \begin{cases}
\makebox[\temp][l]{$- \left( \frac{1}{2} x_6^{2}+ x_7 \right) { x_1} + x_3 \left( \sqrt{2} \left( 1-x_6 \right) x_4 - x_5 x_6 \right) + \frac{1}{2} x_3^{2} x_6$}\\[.25ex]\quad
+\left( x_6^{2} \left( \frac{1}{4} x_6^{2} + x_7 \right) - x_6 +x_7^{2} \right)\cdot \left( \frac{1}{2} x_4^{2} - 2 \sqrt{2} x_4 x_5 + 2 x_5^{2} \right)  \,,
	& \mathrm{if\;} \mathfrak{h} = \mathfrak{co}(2) \ltimes \mathfrak{n}\\[.25ex]
2 x_1 x_6 + \frac{1}{2} x_3^{2} + 3 x_3 x_5 + \frac{9}{2} x_5^{2} +  \left(  x_6^{2} \left( x_7 - \frac{64}{9} \right) + \frac{1}{2} x_7 - 3 \right) x_4^{2}  \,,
	& \mathrm{if\;} \mathfrak{h} = \mathfrak{co}(2) \ltimes Z_0\\[.25ex]
\sqrt{2}  a^3 x_4 x_6^2 \left(x_3 - \frac{3}{2} x_5 \right) +  a^2 x_6 \left( x_1 - x_6 \left( \frac{1}{2} x_3^2 + 7 x_4^2\right)\right)\\[.25ex]\quad
	+  a \left( -x_2 x_6 + \sqrt{2} x_4 \left[ x_3 - \frac{3}{2} x_5 + x_6^2 \left(2x_3+ \frac{5}{2} x_5 \right)\right]\right) \\[.25ex]\qquad
	+ \frac{1}{2} \left( x_3^2 + x_5^2 \right) - 2 x_4^2 + \frac{1}{2} x_3^2 x_6^2 - x_3 x_5 \left(1 + x_6^2\right) \\[.25ex]\qquad\quad
	+ \left( \left(\frac{1}{\sqrt{2}}  a x_4 + x_3 \right) x_5 + x_4^2 \right) \cdot \mathrm{e}^{-\frac{1}{2} \left( a^2 + 1 \right) x_6^2}  \,,
	& \mathrm{if\;} \mathfrak{h} = \mathbb{R}\cdot C_a \ltimes \mathfrak{n}\\[.25ex]
-\frac{ a \left( x_1  a- x_2 \right)}{ a^2 + 1} + \frac{\sqrt{2}  a x_4 \left(  a^2 - \frac{1}{3}\right) \left(  x_3 + 3 x_5 \right)}{\left( a^2 + 1 \right)^2} \\[.25ex]\quad
	- \frac{ x_4^2 \left( 114  a^2  +10 \right)}{36 \left(  a^2 + 1 \right)^2} -\frac{9 x_3^2 + 54 x_3 x_5 + 81 x_5^2}{36 \left(  a^2 + 1 \right)}  \,,
	& \mathrm{if\;} \mathfrak{h} = \mathbb{R}\cdot C_a \ltimes Z_0\\[.25ex]
\end{cases}\\[1ex]
%
%
r_7 &= \begin{cases}
\makebox[\temp][l]{$- \left( \frac{1}{2} x_6^{2}+ x_7 \right)  x_1 -  \left( x_6^{2}+2 x_7 \right)  x_2 + \left( \frac{3}{2} x_6 - 1 \right) x_3^{2} + \left( 2 - 3 x_6 \right) x_3 x_5$}\\[.25ex]\quad + \left( x_6^{2} \left( \frac{1}{4} x_6^{2} + x_7 \right) - x_6 +x_7^{2} \right) \cdot \left( \frac{1}{2} x_4^{2} - 2 \sqrt{2} x_3 x_4 \right)\\[.25ex]\qquad
 + \sqrt{2} \left( 1 -x_6 \right) x_3 x_4 + \frac{3}{2} \left( 1 -  x_6 \right) x_4^{2}  \,,
	& \mathrm{if\;} \mathfrak{h} = \mathfrak{co}(2) \ltimes \mathfrak{n}\\[.25ex]
\left( \frac{2}{3} x_6^{2}-1 \right) x_1 x_6^{2} + 4 x_2 x_6 + \frac{1}{6} x_3^{2} x_6^{3} + \left( 2 \sqrt{2} \left(1 - \frac{8}{3} x_6^{2} \right) x_4 + x_5 x_6^{3} \right) x_3\\[.25ex]\quad
+\left[\frac{1}{3} \left( x_7 - \frac{64}{9} \right) x_6^{4} + \left( \frac{13}{3} - \frac{1}{2} x_7 \right) x_6^{2} + \frac{1}{2} \left( 1 - 8 x_7 \right) \right] x_4^{2} x_6 \\[.25ex]\qquad
+\frac{2 \sqrt{2}}{3} \left( 5 - 8 x_6^{2} \right) x_4 x_5 + \frac{3}{2} x_5^{2} x_6^{3}  \,,
	& \mathrm{if\;} \mathfrak{h} = \mathfrak{co}(2) \ltimes Z_0\\[.25ex]
\left[ 2  a^3 x_4^2 x_6^2 +  a^2 \left(x_2 x_6 + \sqrt{2} x_4 x_6^2 \left(4 x_3 - 5x_5 \right)\right)\right. \\[.25ex]\quad
	+ 2  a \left(x_1 x_6 - x_3 x_6^2 \left( 3 x_3 - 5 x_5 \right) + x_4^2 \left( 1 - 2 x_6^2 \right)\right) \\[.25ex]\qquad
	-\left.  x_2 x_6 - \sqrt{2} x_4 x_5 x_6^2 + \sqrt{2} x_4 \left( 2 x_3 - x_5 \right) \right] \mathrm{e}^{\frac{1}{2} \left( a^2 + 1 \right) x_6^2} \\[.25ex] \qquad\quad
	+2\sqrt{2} x_4 x_5 + \frac{1}{2}  a x_5^2  \,,
	& \mathrm{if\;} \mathfrak{h} = \mathbb{R}\cdot C_a \ltimes \mathfrak{n}\\[.25ex]
\left[-\left( 2 x_1 a+ x_2 a^{2}- x_2 \right)\right.\\\quad
\left.+ 3 x_4\frac{   \left(  \sqrt{2} \left(\left( \frac {7}{9} x_5+ x_3 \right)  a^{2}
-\frac{7}{9} x_3- x_5 \right)\right)+ \frac{1}{3} a x_4 \left(  a^{2}-3 \right)   }{ \left(  a^{2}+1 \right) }\right] \frac{\mathrm{e}^{-x_6}}{ a^{2}+1}  \,,
	& \mathrm{if\;} \mathfrak{h} = \mathbb{R}\cdot C_a \ltimes Z_0\\[.25ex]
\end{cases}\\[1ex]
%
%
s_4 &= \begin{cases}
\makebox[\temp][l]{$ \frac{1}{\sqrt{2}} \left( x_6^{2}  +2 x_7\right) \left(x_3-x_5 \right)  \,,$} 
	& \mathrm{if\;} \mathfrak{h} = \mathfrak{co}(2) \ltimes \mathfrak{n}\\[.25ex]
-2 \sqrt{2} \left( x_3 + \frac{1}{3} x_5 \right) x_6 - 2 x_4 x_7  \,,
	& \mathrm{if\;} \mathfrak{h} = \mathfrak{co}(2) \ltimes Z_0\\[.25ex]
\sqrt{2} x_6 \left(x_3 + \sqrt{2}  a x_4 \right) \,,
	& \mathrm{if\;} \mathfrak{h} = \mathbb{R}\cdot C_a \ltimes \mathfrak{n}\\[.25ex]
-\frac{1}{3}\frac{\sqrt{2} \left( 3 x_3+ x_5 \right) }{ a^{2}+1}  \,,
	& \mathrm{if\;} \mathfrak{h} = \mathbb{R}\cdot C_a \ltimes Z_0\\[.25ex]
\end{cases}\\[1ex]
%
%
s_5 &= \begin{cases}
\makebox[\temp][l]{$-\frac{1}{\sqrt{2}} \left(  x_6^{2}+2x_7 \right) x_4 \,,$}
	& \mathrm{if\;} \mathfrak{h} = \mathfrak{co}(2) \ltimes \mathfrak{n}\\[.25ex]
   a x_3 x_6 \,,
	& \mathrm{if\;} \mathfrak{h} = \mathbb{R}\cdot C_a \ltimes \mathfrak{n}\\[.25ex]
0	 \,,
	& \mathrm{if\;} \mathfrak{h} = \mathfrak{a} \ltimes Z_0\\[.25ex]
\end{cases}\\[1ex]
%
%
s_7 &= \begin{cases}
\makebox[\temp][l]{$ \left( \frac{1}{2} x_6^{2}+ x_7 \right)  x_1- \frac{5}{2} \left( x_6^{2} \left( \frac{1}{4} x_6^{2} + x_7 \right) - x_6 +x_7^{2} \right) x_4^{2} - \frac{3}{2} x_4^{2} x_6 \,,$}
	& \mathrm{if\;} \mathfrak{h} = \mathfrak{co}(2) \ltimes \mathfrak{n}\\[.25ex]
-2 x_1 x_6  - x_2 x_6^{2} + \frac{9}{2} x_3^{2} + \left( 3 x_5 + \frac{8 \sqrt{2}}{3} x_4 x_6^{3} \right) x_3 + \frac{1}{2} x_5^{2} + \frac{8\sqrt{2}}{9} x_4 x_5 x_6^{3}  \,,
	& \mathrm{if\;} \mathfrak{h} = \mathfrak{co}(2) \ltimes Z_0\\[.25ex]
\frac{1}{2}\left(x_3^2 + x_5^2 \right) 
	+ \left[ \left( x_1 +  a x_2 \right) x_6 + \left( \left( 5  a^2 - 3\right) x_6^2 + 1 \right) x_4^2 \right] \cdot \mathrm{e}^{\frac{1}{2} \left( a^2 + 1 \right) x_6^2} \,,
 	& \mathrm{if\;} \mathfrak{h} = \mathbb{R}\cdot C_a \ltimes \mathfrak{n}\\[.25ex]
-\left( x_1 +a x_2 + \frac{9}{4} x_3^{2} +\frac{3}{2} x_3 x_5 + \frac{1}{4} x_5^{2} - \frac{5}{2} x_4^{2} \right) \frac{\mathrm{e}^{-x_6}}{a^{2} + 1}\\[.25ex]\quad
+ \left( 4 \sqrt{2} a \left(x_3 + \frac{1}{3} x_5 \right) x_4 - \frac{14}{3} x_4^{2} \right) \frac{\mathrm{e}^{-x_6}}{\left( a^{2} + 2 \right)^{2}} \,,
	& \mathrm{if\;} \mathfrak{h} = \mathbb{R}\cdot C_a \ltimes Z_0\\[.25ex]
\end{cases}\\[1ex]
\end{align*}
%
%
\settowidth{\temp}{$\quad
	+\left. x_3^2 x_7 \cdot \mathrm{e}^{-x_6 x_7} - 4 x_6 \left( x_7 - 1 \right) \left( x_6 x_7 \left( x_7 - 1 \right) + 3 x_7 - 2 \right) \cdot \mathrm{e}^{x_6 x_7} \right] \,,$}%
\begin{align*}
t_6 &= \begin{cases}
	\makebox[\temp][l]{$-\frac {4\sqrt{2} a x_4}{ a^{2}+1} \,,$}	& \mathrm{if\;} \mathfrak{h} = \mathbb{R}\cdot C_a \ltimes Z_0\\[.25ex]
0	& \mathrm{else\;}\\[.25ex]
\end{cases}\\[1ex]
%
%
t_7 &= \begin{cases}
\makebox[\temp][l]{$\frac{1}{\sqrt{2}}  \left(  x_6^{2}+2x_7 \right)  x_4 \,,$}
	& \mathrm{if\;} \mathfrak{h} = \mathfrak{co}(2) \ltimes \mathfrak{n}\\[.25ex]
-\frac{8\sqrt{2}}{3} x_4 x_6 \,,
	& \mathrm{if\;} \mathfrak{h} = \mathfrak{co}(2) \ltimes Z_0\\[.25ex]
- \left(  a \left( x_3 - x_5 \right) - \sqrt{2} x_4 \right) x_6 \cdot \mathrm{e}^{\frac{1}{2} \left( a^2 + 1 \right) x_6^2} \,,
	& \mathrm{if\;} \mathfrak{h} = \mathbb{R}\cdot C_a \ltimes \mathfrak{n}\\[.25ex]
-\frac{4}{3}\frac{\sqrt{2}x_4 }{ a^{2}+1}\cdot \mathrm{e}^{- x_6} \,,
 	& \mathrm{if\;} \mathfrak{h} = \mathbb{R}\cdot C_a \ltimes Z_0\\[.25ex]
\end{cases}\\[1ex]
%
%
u_6 &= \begin{cases}
\makebox[\temp][l]{$0 \,,$}
	& \mathrm{if\;} \mathfrak{h} = \mathfrak{a} \ltimes \mathfrak{n}\\[.25ex]
\frac{4}{3} x_4 x_6	 \,,
	& \mathrm{if\;} \mathfrak{h} = \mathfrak{co}(2) \ltimes Z_0\\[.25ex]
\frac {3\sqrt{2} a \left( { x_3}-{ x_5} \right) +2 x_4}{3 a^{2}+3} \,,
	& \mathrm{if\;} \mathfrak{h} = \mathbb{R}\cdot C_a \ltimes Z_0\\[.25ex]
\end{cases}\\[1ex]
%
%
u_7 &= \begin{cases}
\makebox[\temp][l]{$-\frac{1}{\sqrt{2}} \left( x_6^{2} + 2 x_7\right) \left(x_3-x_5\right) \,,$}
	& \mathrm{if\;} \mathfrak{h} = \mathfrak{co}(2) \ltimes \mathfrak{n}\\[.25ex]
2 \sqrt{2} x_6 \left(x_3 - x_5 \right) + \frac{4}{9} x_4 x_6^4  \,,
	& \mathrm{if\;} \mathfrak{h} = \mathfrak{co}(2) \ltimes Z_0\\[.25ex]
- \sqrt{2} x_6 \left({ x_3}-{ x_5}+\sqrt{2}  a { x_4} \right) \cdot  \mathrm{e}^{\frac{1}{2}\left(  a^{2}+1 \right)  x_6^{2}}  \,,
	& \mathrm{if\;} \mathfrak{h} = \mathbb{R}\cdot C_a \ltimes \mathfrak{n}\\[.25ex]
\frac {\sqrt{2}\left( { x_3}-{ x_5} \right) +2\,{ x_4}\, a  }{ a^{2}+1}\cdot \mathrm{e}^{- x_6}  \,,
	& \mathrm{if\;} \mathfrak{h} = \mathbb{R}\cdot C_a \ltimes Z_0\\[.25ex]
\end{cases}\\[1ex]
%
%
w_7 &= \begin{cases}
\makebox[\temp][l]{$1 \,,$}	& \mathrm{if\;} \mathfrak{h} = \mathfrak{co}(2) \ltimes \mathfrak{n}\\[.25ex]
\frac{x_6^3}{3}  \,,
	& \mathrm{if\;} \mathfrak{h} = \mathfrak{co}(2) \ltimes Z_0\\[.25ex]
0	 \,,
	& \mathrm{if\;} \mathfrak{a} = \mathbb{R}\cdot C_a \\[.25ex]
\end{cases}\\[1ex]
z_7 &= \begin{cases}
\makebox[\temp][l]{$0 \,,$}	& \mathrm{if\;} \mathfrak{a} = \mathfrak{co}(2)\\[.25ex]
\frac{1}{2} \left(  a^{2}+1 \right)  x_6^{2}	\,,
	& \mathrm{if\;} \mathfrak{h} = \mathbb{R}\cdot C_a \ltimes \mathfrak{n}\\[.25ex]
-x_6 	 \,,
	& \mathrm{if\;} \mathfrak{h} = \mathbb{R}\cdot C_a \ltimes Z_0\\[.25ex]
\end{cases}\\[1ex]
\end{align*}
\normalsize

%% file: Table_Type_II_3.tex
\begin{center}
\textsc{Type II 3}
\end{center}
$\mathfrak{a} = \mathfrak{d}$ and $\mathfrak{h} = \mathfrak{a}\ltimes \mathfrak{n_1}$, where $\mathfrak{n}_1 \in \left\lbrace \mathfrak{n},\mathfrak{n}(1,3),\mathfrak{n}(2,3),\mathfrak{n}(1,2,3),\mathfrak{n}(1,2,4) \right\rbrace$.
\small
\settowidth{\temp}{$\quad
	+\left. x_3^2 x_7 \cdot \mathrm{e}^{-x_6 x_7} - 4 x_6 \left( x_7 - 1 \right) \left( x_6 x_7 \left( x_7 - 1 \right) + 3 x_7 - 2 \right) \cdot \mathrm{e}^{x_6 x_7} \right] \,, $}
%
%
\begin{align*}
r_5 &= \begin{cases}
\makebox[\temp][l]{$-2 x_3 x_7 \,,$}	& \mathrm{if\;} \mathfrak{n}_1 = \mathfrak{n}\\[.25ex]
0 \,,
	& \mathrm{if\;} \mathfrak{n}_1 = \mathfrak{n}(1,3)\\[.25ex]
-\frac{1}{2} x_3 x_7 \left(2 + x_7\right) \,,
	& \mathrm{if\;} \mathfrak{n}_1 = \mathfrak{n}(2,3)\\[.25ex]
1 - x_3 x_7 \,,
	& \mathrm{if\;} \mathfrak{n}_1 = \mathfrak{n}(1,2,3)\\[.25ex]
2 x_3 x_5 \,,
	& \mathrm{if\;} \mathfrak{n}_1 = \mathfrak{n}(1,2,4)\\[.25ex]
\end{cases}\\[1ex]
%
%
r_6 &= \begin{cases}
\makebox[\temp][l]{$- x_1 x_7 - 2 x_3 x_5 x_7^{2} - 2 x_4^{2} x_7^{2} + \frac{1}{2} x_5^{2} + \sqrt{2} x_4 \left( x_3 - \frac{1}{2} x_3^{2} x_7 \right) \,, $}
	& \mathrm{if\;} \mathfrak{n}_1 = \mathfrak{n}\\[.25ex]
- x_1 x_7 + \sqrt{2} x_4 \left( x_6 x_7^2 + 1\right) \cdot \mathrm{e}^{x_6 x_7} \,,
	& \mathrm{if\;} \mathfrak{n}_1 = \mathfrak{n}(1,3)\\[.25ex]
-\frac{1}{2} x_1 x_7 \left( 2 + x_7 \right) + \sqrt{2} x_3 x_4 \,,
	& \mathrm{if\;} \mathfrak{n}_1 = \mathfrak{n}(2,3)\\[.25ex]
-x_1 x_7 +\frac{1}{2} \left( x_3^{2}\left(1+ \sqrt{2}  x_4 x_7\right)- \mathrm{e}^{- x_3}\right)  \,,
	& \mathrm{if\;} \mathfrak{n}_1 = \mathfrak{n}(1,2,3)\\[.25ex]
	 x_1 x_5 \,,
	& \mathrm{if\;} \mathfrak{n}_1 = \mathfrak{n}(1,2,4)\\[.25ex]
\end{cases}\\[1ex]
%
%
r_7 &= \begin{cases}
\makebox[\temp][l]{$2 x_4^2 - 2 x_4 x_7 \left( x_3 x_4 + 2 \sqrt{2} x_5 x_7 \right) \,,$}
	& \mathrm{if\;} \mathfrak{n}_1 = \mathfrak{n}\\[.25ex]
\frac{1}{2} x_5 \left[ 2 x_3  \left( x_6 x_7 \left(x_7 - 1\right) - 1 \right) \right. \\[.25ex]\quad
	+\left. x_3^2 x_7 \cdot \mathrm{e}^{-x_6 x_7} - 4 x_6 \left( x_7 - 1 \right) \left( x_6 x_7 \left( x_7 - 1 \right) + 3 x_7 - 2 \right) \cdot \mathrm{e}^{x_6 x_7} \right] \,,
	& \mathrm{if\;} \mathfrak{n}_1 = \mathfrak{n}(1,3)\\[.25ex]
2 x_4^2 - x_3 x_5 x_7 \left( 1 + \frac{1}{2} x_6 x_7^2 + x_6 x_7\right) \,,
	& \mathrm{if\;} \mathfrak{n}_1 = \mathfrak{n}(2,3)\\[.25ex]
2 x_3^2 x_5 x_7 + 2 x_3\left( x_4^{2} x_7+\sqrt{2}  x_4- x_5\right) + \sqrt{2} x_4 \cdot \mathrm{e}^{-x_3} \,,
	& \mathrm{if\;} \mathfrak{n}_1 = \mathfrak{n}(1,2,3)\\[.25ex]
2 x_3 x_6 +\frac{1}{2} x_3^{2} x_5^{2} \,,
	& \mathrm{if\;} \mathfrak{n}_1 = \mathfrak{n}(1,2,4)\\[.25ex]
\end{cases}\\[1ex]
\end{align*}
%
%
\settowidth{\temp}{$x_3 x_6 (x_7 -1)+\frac{x_3^{2}}{2} \mathrm{e}^{-x_6 x_7 }  + x_6 \left( 1 - 2 x_6 \left( x_7 -1 \right)^{2}\right) \mathrm{e}^{x_6 x_7} \,, $}
\begin{align*}
s_5 &= \begin{cases}
\makebox[\temp][l]{$\frac{1}{2} x_3^2  - \sqrt{2} x_4 x_7 \,,$}	
	& \mathrm{if\;} \mathfrak{n}_1 = \mathfrak{n}\\[.25ex]
 x_3 x_6 (x_7 -1)+\frac{x_3^{2}}{2} \mathrm{e}^{-x_6 x_7 }  + x_6 \left( 1 - 2 x_6 \left( x_7 -1 \right)^{2}\right) \mathrm{e}^{x_6 x_7} \,,
	& \mathrm{if\;} \mathfrak{n}_1 = \mathfrak{n}(1,3)\\[.25ex]
-\frac{1}{\sqrt{2}} x_4 x_7 \left( 2 + x_7\right) \,,
	& \mathrm{if\;} \mathfrak{n}_1 = \mathfrak{n}(2,3)\\[.25ex]
\frac{1}{2} x_3^{2}-\sqrt{2}\, x_4 x_7 \,,
	& \mathrm{if\;} \mathfrak{n}_1 = \mathfrak{n}(1,2,3)\\[.25ex]
-\frac{1}{2} x_3^{2}+ \sqrt {2}\,x_4 x_5 \,,
	& \mathrm{if\;} \mathfrak{n}_1 = \mathfrak{n}(1,2,4)\\[.25ex]
\end{cases}\\[1ex]
%
%
s_7 &= \begin{cases}
\makebox[\temp][l]{$- x_2 x_3 - x_3^{3} x_5 - \frac{\sqrt{2}}{3} x_4^{3} x_7 \,,$}
	& \mathrm{if\;} \mathfrak{n}_1 = \mathfrak{n}\\[.25ex]
-x_2 \left(x_6 \left(x_7 - 1 \right) + x_3 \cdot \mathrm{e}^{-x_6 x_7} \right) - x_3^3 x_5 \cdot \mathrm{e}^{-2 x_6 x_7}\\[.25ex]\quad
	 +x_5 x_6 \left(\left( 2 x_6 \left( x_7 -2\right) \left( x_7 - 1 \right)^{2} - 5 x_7 + 6 \right) x_6 \cdot  \mathrm{e}^{x_6 x_7} \right.\\[.25ex]\quad +\left.  x_3^{2} \left( \frac{5}{2} - 3 x_7 \right) \cdot \mathrm{e}^{-x_6 x_7} \right) \,,
 	& \mathrm{if\;} \mathfrak{n}_1 = \mathfrak{n}(1,3)\\[.25ex]
 x_2 x_6 x_7 - \frac{1}{\sqrt{2}} x_4 x_5 x_6 x_7^2 \left( 2 +x_7 \right) - \sqrt{2} x_4 x_5 x_7 \\[.25ex]\quad
 	+ x_5^2 x_7^2 \left(\frac{1}{4} x_7^2 + x_7 + 1 \right) \,,
	& \mathrm{if\;} \mathfrak{n}_1 = \mathfrak{n}(2,3)\\[.25ex]
-x_2 x_3 + x_3 x_5 \left( 1 -x_3^2 +2\sqrt{2} x_4 x_7 \right)+ x_5^2 x_7^2 \\[.25ex]\quad
	+ x_4 \left( x_4 + \frac{\sqrt{2}}{3} x_4^2 x_7 - \sqrt{2} x_5 \right) + \frac{1}{2}\left(x_5 -x_4^2\right)\cdot \mathrm{e}^{-x3} \,,
	& \mathrm{if\;} \mathfrak{n}_1 = \mathfrak{n}(1,2,3)\\[.25ex]
x_2 x_3+\sqrt{2}\, x_4 x_6+\frac{1}{\sqrt{2}} x_3 x_4 x_5^{2} - x_3^{3} x_5+\frac{1}{4} x_5^{4} \,,
	& \mathrm{if\;} \mathfrak{n}_1 = \mathfrak{n}(1,2,4)
\end{cases}\\[1ex]
%
%
t_6 &= \begin{cases}
\makebox[\temp][l]{$0 \,,$}
	& \mathrm{if\;} \mathfrak{n}_1 \in \left\lbrace \mathfrak{n},\mathfrak{n}(1,2,4)\right\rbrace\\[.25ex]
 x_7 \left( \mathrm{e}^{ x_6 x_7} - x_3\right) \,,
	& \mathrm{if\;} \mathfrak{n}_1 = \mathfrak{n}(1,3)\\[.25ex]
-\frac{1}{2} x_3 x_7 \left(2 + x_7\right) \,,
	& \mathrm{if\;} \mathfrak{n}_1 = \mathfrak{n}(2,3)\\[.25ex]
1- x_3 x_7	\,,
	& \mathrm{if\;} \mathfrak{n}_1 = \mathfrak{n}(1,2,3)\\[.25ex]
\end{cases}\\[1ex]
%
%
t_7 &= \begin{cases}
\makebox[\temp][l]{$\frac{1}{2} x_3^{2}\,,$}	
	& \mathrm{if\;} \mathfrak{n}_1 \in \left\lbrace \mathfrak{n}, \mathfrak{n}(1,2,3)\right\rbrace\\[.25ex]
 x_3 x_6 (x_7 -1)+\frac{x_3^{2}}{2} \mathrm{e}^{-x_6 x_7 } \,,
	&\hspace{0cm} \mathrm{if\;} \mathfrak{n}_1 = \mathfrak{n}(1,3)\\[.25ex]
- x_3 x_6 x_7	\,,
	&\hspace{0cm} \mathrm{if\;} \mathfrak{n}_1 = \mathfrak{n}(2,3)\\[.25ex]
-\frac{1}{2} x_3^{2}	\,,
	&\hspace{0cm} \mathrm{if\;} \mathfrak{n}_1 = \mathfrak{n}(1,2,4)\\[.25ex]
\end{cases}\\[1ex]
%
%
u_7 &= \begin{cases}
\makebox[\temp][l]{$\sqrt{2}\, x_5 x_7\,,$}
	& \mathrm{if\;} \mathfrak{n}_1 \in \left\lbrace \mathfrak{n}, \mathfrak{n}(1,2,3)\right\rbrace\\[.25ex]
0 \,,
	&\hspace{0cm} \mathrm{if\;} \mathfrak{n}_1 = \mathfrak{n}(1,3)\\[.25ex]
\frac{1}{\sqrt{2}} x_5 x_7 \left(2 + x_7\right)	\,,
	&\hspace{0cm} \mathrm{if\;} \mathfrak{n}_1 = \mathfrak{n}(2,3)\\[.25ex]
-\frac{ x_5^{2}}{\sqrt{2}} \,,
	& \mathrm{if\;} \mathfrak{n}_1 = \mathfrak{n}(1,2,4)
\end{cases}\\[1ex]
%
%
v_6 &= \begin{cases}
\makebox[\temp][l]{$ x_5 x_7 \,,$}	
	& \mathrm{if\;} \mathfrak{n}_1 = \mathfrak{n}(1,3)\\[.25ex]
0
	& \mathrm{else\;} 
\end{cases}\\[1ex]
%
%
v_7 &= \begin{cases}
\makebox[\temp][l]{$x_5 x_6 x_7 \,,$}	& \mathrm{if\;} \mathfrak{n}_1 = \mathfrak{n}(2,3)\\[.25ex]
0
	& \mathrm{else\;}
\end{cases}\\[1ex]
\end{align*}
\normalsize

%% file: Table_Type_II_4.tex
\clearpage
\begin{center}
\textsc{Type II 4}
\end{center}
Let $\mathfrak{a}=\mathbb{R}\cdot \mathrm{diag}(1,\mu)$.\\
%
%
\noindent\textbf{a)}
$\mathfrak{h}=\mathfrak{a}\ltimes \mathfrak{n_1},  \mu \in [-1,1),\gamma=1-\mu$,\\
where $\mathfrak{n}_1 \in \left\lbrace \mathfrak{n},\mathfrak{n}(2,3),\mathfrak{n}(1,2,3),\mathfrak{n}(1,2,4),\mathfrak{n}(1,3,4),\mathfrak{n}(2,3,4)\right\rbrace
$.
\small
\settowidth{\temp}{$\mu x_2 x_6+\frac{1}{2} x_3^2 + \sqrt {2}\gamma x_4 x_5 \cdot \mathrm{e}^{-x_6 x_7} + \left(\mu\gamma - \gamma^2\right) x_3 x_5 x_6^2 \,,$}
\begin{align*}
r_6 &= \begin{cases}
\makebox[\temp][l]{$\sqrt{2}\, x_3 x_4 \,, $}
	& \mathrm{if\;} \mathfrak{n}_1 = \mathfrak{n}\\[.25ex]
\sqrt{2} \gamma x_3 x_4	\,,
	& \mathrm{if\;} \mathfrak{n}_1 \in \left\lbrace \mathfrak{n}(2,3),\mathfrak{n}(1,2,3) \right\rbrace\\[.25ex]
\frac{1}{2} x_5^{2} \,,
	& \mathrm{if\;} \mathfrak{n}_1 = \mathfrak{n}(1,2,4)\\[.25ex]
 \sqrt{2}\,x_4 x_5 \cdot \mathrm{e}^{x_6 x_7}\,,
	& \mathrm{if\;} \mathfrak{n}_1 = \mathfrak{n}(1,3,4)\\[.25ex]
\sqrt{2}\,\gamma x_3 x_4+\frac{1}{2} x_5^{2} \,,
	& \mathrm{if\;} \mathfrak{n}_1 = \mathfrak{n}(2,3,4)\\[.25ex]
\end{cases}\\[1ex]
%
%
r_7 &=  \begin{cases}
\makebox[\temp][l]{$ x_1 x_6 +2 \left(x_3 x_5+ x_4^{2}\right)\cdot\mathrm{e}^{-x_6 x_7} \,,$}
	& \mathrm{if\;} \mathfrak{n}_1 = \mathfrak{n}\\[.25ex]
 x_1 x_6 +\gamma \left(  x_3 x_5+2 x_4^{2} \right) \cdot \mathrm{e}^{-x_6 x_7} \,,
	& \mathrm{if\;} \mathfrak{n}_1 \in \{ \mathfrak{n}(2,3),\mathfrak{n}(2,3,4)\}\\[.25ex]
 x_1 x_6 +2\gamma\left( x_3 x_5+ x_4^{2}\right)\cdot\mathrm{e}^{-x_6 x_7} \,,
	& \mathrm{if\;} \mathfrak{n}_1 = \mathfrak{n}(1,2,3)\\[.25ex]
 x_1 x_6 + x_3^{2} \,,
	& \mathrm{if\;} \mathfrak{n}_1 = \mathfrak{n}(1,2,4)\\[.25ex]
 x_1 x_6 + x_5^{2} \,,
	& \mathrm{if\;} \mathfrak{n}_1 = \mathfrak{n}(1,3,4)\\[.25ex]
\end{cases}\\[1ex]
%
%
s_5 &= \begin{cases}
0	\,,
	& \mathrm{if\;} \mathfrak{n}_1 \in \left\lbrace \mathfrak{n}(2,3),\mathfrak{n}(2,3,4)\right\rbrace\\[.25ex]
\makebox[\temp][l]{$ \gamma x_3 x_6 \,,	$}
	& \mathrm{else\;}
\end{cases}\\[1ex]
%
%
s_7 &= \begin{cases}
\makebox[\temp][l]{$\mu x_2 x_6 + \sqrt{2} x_4 x_5 \cdot \mathrm{e}^{-x_6 x_7} + \left(2\mu - 1\right)\gamma x_3 x_5 x_6^2 \,,$}
 	& \mathrm{if\;} \mathfrak{n}_1 = \mathfrak{n}\\[.25ex]
\mu x_2 x_6+\sqrt{2}\gamma x_4 x_5 \cdot \mathrm{e}^{-x_6 x_7} \,,
	& \mathrm{if\;} \mathfrak{n}_1 \in \{ \mathfrak{n}(2,3),\mathfrak{n}(2,3,4)\} \\[.25ex]
\mu x_2 x_6+\frac{1}{2} x_3^2 + \sqrt {2}\gamma x_4 x_5 \cdot \mathrm{e}^{-x_6 x_7} + \left(\mu\gamma - \gamma^2\right) x_3 x_5 x_6^2 \,,
	& \mathrm{if\;} \mathfrak{n}_1 =  \mathfrak{n}(1,2,3) \\[.25ex]
\mu x_2 x_6+ \sqrt {2} x_3 x_4 + \left(\mu\gamma - \gamma^2\right) x_3 x_5 x_6^{2}  \,,
 & \mathrm{if\;} \mathfrak{n}_1 = \mathfrak{n}(1,2,4)\\[.25ex]
\mu x_2 x_6 +\frac{1}{2} x_3^{2}+  \left( \mu\gamma-\gamma^2 \right) x_3 x_5 x_6^{2} \,,
	& \mathrm{if\;} \mathfrak{n}_1 = \mathfrak{n}(1,3,4)\\[.25ex]
\end{cases}\\[1ex]
%
%
v_7 &= \begin{cases}
\makebox[\temp][l]{$-\gamma x_5 x_6 \,,$}
 	& \mathrm{if\;} \mathfrak{n}_1 \in \{ \mathfrak{n}(2,3),\mathfrak{n}(2,3,4)\}\\[.25ex]
0
	& \mathrm{else\;} \\[.25ex]
\end{cases}\\[1ex]
%
%
t_7 &= \makebox[\temp][l]{$\gamma x_3 x_6, \qquad u_7 = 0,  \qquad w_6 = x_6 x_7 $} \qquad \mathrm{in\;all\;cases} \\
\end{align*}
%
%
\textbf{b)}
$\mathfrak{h}=h\left(\mathrm{diag}(1,\frac{1}{2}),(1,0,0,0),0\right)\ltimes \mathfrak{n_1}$, where $\mathfrak{n}_1 \in \left\lbrace \mathfrak{n}(2,3),\mathfrak{n}(2,3,4)\right\rbrace$.
\settowidth{\temp}{$\frac{1}{2}\,{ x_3}\,{ x_6} \qquad u_7 = 0 \qquad v_7 = -\frac{1}{2}\,{ x_5}\,{ x_6} \qquad w_6=\mathrm{e}^{x_6 x_7}\qquad$}
%
%
\begin{align*}
r_6 &= \begin{cases}
\makebox[\temp][l]{$\frac{1}{\sqrt {2}}x_3 x_4 \,,$}	
	& \mathrm{if\;} \mathfrak{n}_1 = \mathfrak{n}(2,3)\\[.25ex]
\frac{1}{\sqrt {2}}x_3 x_4+\frac{1}{2} x_5^{2} \,,
	&	 \mathrm{if\;} \mathfrak{n}_1 =  \mathfrak{n}(2,3,4)
\end{cases}\\
%
%
r_7 &= 
 \makebox[\temp][l]{$ x_1 x_6+ \left( \frac{1}{2} x_3 x_5+ x_4^{2}\right)\cdot  \mathrm{e}^{-x_6 x_7}$} 
	\qquad \mathrm{in\;both\;cases} \\
%
%
s_7 &= 
\makebox[\temp][l]{$\left( \frac{1}{2}{ x_2}+{ x_3} \right) { x_6}+\frac{1}{\sqrt {2}} x_4 x_5 \cdot  \mathrm{e}^{-x_6 x_7}$}
	\qquad \mathrm{in\;both\;cases} \\
t_7 &= 
\makebox[\temp][l]{$ \frac{1}{2}\,{ x_3}\,{ x_6}, \qquad u_7 = 0, \qquad v_7 = -\frac{1}{2}\,{ x_5}\,{ x_6}, \qquad w_6=x_6 x_7 $}	\qquad \mathrm{in\;both\;cases} 
\end{align*}
\newpage
%
%
\noindent\textbf{c)}
$\mu=0$ and $\mathfrak{h}=\mathfrak{a}\ltimes \mathfrak{n}(2,4)$ or $\mathfrak{h}=\mathbb{R}\cdot h\left(\mathrm{diag}(1,0),(0,1,0,0),0\right) \ltimes \mathfrak{n}_1$,  $\mathfrak{n}_1 \in \left\lbrace \mathfrak{n}(1,4), \mathfrak{n}(3,4),\mathfrak{n}(1,3,4)\right\rbrace$.
\settowidth{\temp}{$ x_1 x_5+ x_3 \left( 1-\frac{1}{2} x_5^{2} \right)  x_5+ x_4^{2} \,,$}
%
%
\begin{align*}
r_5 &= \begin{cases}
\makebox[\temp][l]{$2 x_3 x_5 - \frac{1}{3} x_5 x_7^3 \,,$}
	& \mathrm{if\;} \mathfrak{h} = \mathfrak{a}\ltimes \mathfrak{n}(2,4)\\[.25ex]
- x_3 x_7 \,,
	& \mathrm{if\;} \mathfrak{n}_1 = \mathfrak{n}(1,4)\\[.25ex]
 x_3 x_5 \,,
	& \mathrm{if\;} \mathfrak{n}_1 \in\left\lbrace  \mathfrak{n}(3,4),\mathfrak{n}(1,3,4)\right\rbrace
\end{cases}\\[1ex]
r_6 &= \begin{cases}
\makebox[\temp][l]{$x_1 x_5- x_3 x_7 \,,$}	
	& \mathrm{if\;} \mathfrak{h} = \mathfrak{a}\ltimes \mathfrak{n}(2,4)\\[.25ex]
-   \left(x_1 +x_3\right)x_7- x_3 x_5 x_7^{2}	\,,
	& \mathrm{if\;} \mathfrak{n}_1 = \mathfrak{n}(1,4)\\[.25ex]
 x_1 x_5+ x_3 \left( 1-\frac{1}{2} x_5^{2} \right)  x_5+ x_4^{2} \,,
	& \mathrm{if\;} \mathfrak{n}_1 \in\left\lbrace  \mathfrak{n}(3,4),\mathfrak{n}(1,3,4)\right\rbrace
\end{cases}\\[1ex]
%
%
r_7 &= \begin{cases}
\makebox[\temp][l]{$-\sqrt {2} x_4 x_7 \,,$}
	& \mathrm{if\;} \mathfrak{h} = \mathfrak{a}\ltimes \mathfrak{n}(2,4)\\[.25ex]
- x_3 x_5-\sqrt{2} x_4 x_7	\,,
	& \mathrm{if\;} \mathfrak{n}_1 = \mathfrak{n}(1,4)\\[.25ex]
2\sqrt {2}\, x_4 x_5	 \,,
	& \mathrm{if\;} \mathfrak{n}_1 \in\left\lbrace  \mathfrak{n}(3,4),\mathfrak{n}(1,3,4)\right\rbrace
\end{cases}\\[1ex]
%
%
s_5 &= \begin{cases}
\makebox[\temp][l]{$\sqrt{2} x_4 x_5 - x_5^{3} x_7 + x_5 x_6 x_7^2 \,,$}
	& \mathrm{if\;} \mathfrak{h} = \mathfrak{a}\ltimes \mathfrak{n}(2,4)\\[.25ex]
0 \,,
	& \mathrm{if\;} \mathfrak{n}_1 \in \left\lbrace \mathfrak{n}(1,4),\mathfrak{n}(3,4),\mathfrak{n}(1,3,4)\right\rbrace
\end{cases}\\[1ex]
%
%
s_7 &= \begin{cases}
\makebox[\temp][l]{$0 \,,$}	
	& \mathrm{if\;} \mathfrak{h} = \mathfrak{a}\ltimes \mathfrak{n}(2,4) \;\mathrm{or}\;\mathfrak{n}_1 = \mathfrak{n}(1,4)\\[.25ex]
\frac{1}{2} x_5^{2}+\frac{1}{2} x_6^{2} \,,
	& \mathrm{if\;} \mathfrak{n}_1 = \mathfrak{n}(3,4)\\[.25ex]
\frac{1}{2}  x_5^{2}+ x_3 x_6 \,,
	& \mathrm{if\;} \mathfrak{n}_1 = \mathfrak{n}(1,3,4)
\end{cases}\\[1ex]
%
%
 t_7 &= \makebox[\temp][l]{$0, \qquad u_6 = 0$}\qquad \mathrm{in\;all\;cases}\\[1ex]
%
%
u_7 &= \begin{cases}
\makebox[\temp][l]{$-\frac{1}{\sqrt {2}} x_5^{2} + \sqrt{2} x_6 x_7 \,,$}	
	& \mathrm{if\;} \mathfrak{h} = \mathfrak{a}\ltimes \mathfrak{n}(2,4)\\[.25ex]
\sqrt{2} x_6 x_7 \,,
	& \mathrm{if\;} \mathfrak{n}_1 = \mathfrak{n}(1,4)\\[.25ex]
0 \,,
	& \mathrm{if\;} \mathfrak{n}_1 \in\left\lbrace  \mathfrak{n}(3,4),\mathfrak{n}(1,3,4)\right\rbrace
\end{cases}\\[1ex]
%
%
v_6 &= \begin{cases}
\makebox[\temp][l]{$0 \,,$}
	& \mathrm{if\;} \mathfrak{h} = \mathfrak{a}\ltimes \mathfrak{n}(2,4)\\[.25ex]
 x_5  x_7 \,,
	& \mathrm{if\;} \mathfrak{n}_1 = \mathfrak{n}(1,4)\\[.25ex]
-\frac{1}{2}  x_5^{2} \,,
	& \mathrm{if\;} \mathfrak{n}_1 \in\left\lbrace  \mathfrak{n}(3,4),\mathfrak{n}(1,3,4)\right\rbrace
\end{cases}\\[1ex]
\end{align*}
\normalsize

%% file: Table_Type_II_5.tex
\begin{center}
\textsc{Type II 5}
\end{center}
%
%
\textbf{a)}
Let $\mathfrak{h} = \mathfrak{n}_1$, where $\mathfrak{n}_1 \in \left\lbrace \mathfrak{n},\mathfrak{n}(1,3),\mathfrak{n}(2,3),\mathfrak{n}(1,2,4),\mathfrak{n}(2,3,4),Z_1,Z_2,Z_3,Z_4,Z_5 \right\rbrace
$.
\small
\settowidth{\temp}{$\left(\alpha{ x_3}+\sqrt{2}\;{ x_4}+{\frac {{ x_5}}{\alpha}} \right)\cdot  \left({x_7} +{{x_6}}^{2}-\alpha{x_6}({x_7}+1) \right)$}
%
%
\begin{align*}
r_6 &= \begin{cases}
\makebox[\temp][l]{$\sqrt{2} x_3 \left(  x_4+ x_6 \right)$} 
	& \mathrm{if\;} \mathfrak{n}_1 = \mathfrak{n}\\[.25ex]
-\sqrt{2} x_4 x_7
	& \mathrm{if\;} \mathfrak{n}_1 = \mathfrak{n}(1,3)\\[.25ex]
	 x_3 \left( \frac{1}{2} x_6^{2}-x_7 \right) 
	& \mathrm{if\;} \mathfrak{n}_1 = \mathfrak{n}(2,3)\\[.25ex]
	\sqrt{2} x_4 x_5 
	& \mathrm{if\;} \mathfrak{n}_1 = \mathfrak{n}(1,3,4)\\[.25ex]
	-x_3 x_7 + \sqrt{2} x_4 x_5
	& \mathrm{if\;} \mathfrak{n}_1 = \mathfrak{n}(2,3,4)\\[.25ex]
	\frac{1}{2} x_5^{2} 
	& \mathrm{if\;} \mathfrak{n}_1 = Z_1 \\[.25ex]
	-\sqrt{2} x_4 x_6 x_7 - \left( x_3 - x_5 \right) x_7
 	& \mathrm{if\;} \mathfrak{n}_1 = Z_2\\
 -\frac{1}{2}x_7^2 \left( \sqrt{2} \alpha x_4 + x_5 \right)
	& \mathrm{if\;} \mathfrak{n}_1 = Z_3\\
2 x_4 x_6 x_7
	& \mathrm{if\;} \mathfrak{n}_1 = Z_4\\
\sqrt{2} \kappa x_4 x_6 + \frac{1}{2} x_5^{2} 
	& \mathrm{if\;} \mathfrak{n}_1 = Z_5
\end{cases}\\[1ex]
\end{align*}
%
%
\begin{align*}
r_7 &=\begin{cases}
\makebox[\temp][l]{$2\left(x_3 x_5+x_4^2+2 x_4 x_6 \right)	\,,$} 
	& \mathrm{if\;} \mathfrak{n}_1 = \mathfrak{n}\\[.25ex]
- x_5 x_7	\,,
	& \mathrm{if\;} \mathfrak{n}_1 = \mathfrak{n}(1,3)\\[.25ex]
\sqrt{2} x_4 \left(x_6^{2}- x_7 \right) 	\,, 
	& \mathrm{if\;} \mathfrak{n}_1 = \mathfrak{n}(2,3)\\[.25ex]
 x_5^{2}	 \,,
	& \mathrm{if\;} \mathfrak{n}_1 = \mathfrak{n}(1,3,4)\\[.25ex]
 -2\sqrt{2} x_4 x_7\,,
	& \mathrm{if\;} \mathfrak{n}_1 = \mathfrak{n}(2,3,4)\\[.25ex]
\sqrt{2} x_4 x_6 +2 x_5 x_7 \,,
	& \mathrm{if\;} \mathfrak{n}_1 = Z_1 \\[.25ex]
-2\sqrt{2} x_4 x_7 - x_5 x_6 x_7 \,,
	& \mathrm{if\;} \mathfrak{n}_1 = Z_2\\
\frac{\alpha}{2}\left(2 x_3 x_6 - x_5 x_7^2 \right)\,,
	& \mathrm{if\;} \mathfrak{n}_1 = Z_3\\
-\sqrt{2} x_6 \left(\frac{\beta}{s} x_4 + \frac{x_7}{\beta} \left( 2 s x_3 - \frac{\beta^{2} + s}{\beta} x_5 \right)\right)\,,
	& \mathrm{if\;} \mathfrak{n}_1 = Z_4\\
\left( 2 x_3 + \kappa x_5 \right) x_6 + \frac{1}{\sqrt{2}} \kappa x_4 x_6^{2}\,,
	& \mathrm{if\;} \mathfrak{n}_1 = Z_5
\end{cases}\\[1ex]
%
%
s_5 &= \begin{cases}
 \makebox[\temp][l]{$ \frac{1}{2} x_7^2 -2 x_5 x_6 \,,$} 
	& \mathrm{if\;} \mathfrak{n}_1 = \mathfrak{n}(2,3,4)\\
0   & \mathrm{else} 
\end{cases}\\[1ex]
%
%
s_7 &= \begin{cases}
\makebox[\temp][l]{$\sqrt{2} x_5 \left( x_4+ x_6 \right) \,,$} 
	& \mathrm{if\;} \mathfrak{n}_1 = \mathfrak{n}\\[.25ex]
\frac{1}{2} x_3^{2}	\,,
	& \mathrm{if\;} \mathfrak{n}_1 = \mathfrak{n}(1,3)\\[.25ex]
\frac{1}{2} x_5 x_6^{2}	\,,
	& \mathrm{if\;} \mathfrak{n}_1 = \mathfrak{n}(2,3)\\[.25ex]
x_3 x_6	\,,
	& \mathrm{if\;} \mathfrak{n}_1 = \mathfrak{n}(1,3,4)\\[.25ex]
0    \,,
	& \mathrm{if\;} \mathfrak{n}_1 = \mathfrak{n}(2,3,4)\\[.25ex]
\left( x_3+x_5 \right) x_6	\,,
	& \mathrm{if\;} \mathfrak{n}_1 = Z_1 \\[.25ex]
-x_5 x_7   \,,
	& \mathrm{if\;} \mathfrak{n}_1 = Z_2\\
\left( \sqrt{2} \alpha x_4 + x_3 \right) x_6 \,,
	& \mathrm{if\;} \mathfrak{n}_1 = Z_3\\
-x_6 \left( \frac{2 s}{\beta} x_4 x_7 + \frac{1}{s} \left( \beta x_5 - s x_3 \right)\right)  \,,
    & \mathrm{if\;} \mathfrak{n}_1 = Z_4\\
 \frac{1}{2}x_6^{2} \left(x_3 + \kappa x_5 \right) +  \sqrt{2} x_4 x_6 \,,
    & \mathrm{if\;} \mathfrak{n}_1 = Z_5
\end{cases}\\[1ex]
%
%
t_6 &= \begin{cases}
 \makebox[\temp][l]{$ \frac{\alpha}{6} x_7^3 \,,$} & \mathrm{if\;} \mathfrak{n}_1 = Z_3\\
0   & \mathrm{else} 
\end{cases}\\[1ex]
%
%
t_7 &= \begin{cases}
\makebox[\temp][l]{$- x_6 x_7 \,,$} 
	& \mathrm{if\;} \mathfrak{n}_1 = \mathfrak{n}(1,3)\\[.25ex]
-2 x_6 x_7	\,,
	& \mathrm{if\;} \mathfrak{n}_1 = Z_1 \\[.25ex]
-\frac{1}{2} x_6^{2} x_7 \,,
	& \mathrm{if\;} \mathfrak{n}_1 = Z_2\\
\frac{1}{\sqrt{2}}\frac{\left( \beta^2 - s \right)}{\beta}  x_6^{2} x_7 \,,
	& \mathrm{if\;} \mathfrak{n}_1 = Z_4\\
    \frac{\kappa}{2} x_6^{2}	 \,,
	& \mathrm{if\;} \mathfrak{n}_1 = Z_5 \\[.25ex]
0    & \mathrm{else}
\end{cases}\\[1ex]
%
%
u_6 &= \begin{cases}
\makebox[\temp][l]{$\sqrt{2} x_6 \,,$} 
	& \mathrm{if\;} \mathfrak{n}_1 = \mathfrak{n}\\[.25ex]
	-\frac{\kappa}{\sqrt{2}} x_6^{2} x_7	  \,,
	& \mathrm{if\;} \mathfrak{n}_1 = Z_5 \\[.25ex]
0    & \mathrm{else}
\end{cases}\\[1ex]
%
%
u_7 &= \begin{cases}
\makebox[\temp][l]{$\sqrt{2}x_6 x_7 \,,$} 
	& \mathrm{if\;} \mathfrak{n}_1 = \mathfrak{n}(2,3)\\[.25ex]
\frac{1}{\sqrt{2}} x_6^{2}	 \,,
	& \mathrm{if\;} \mathfrak{n}_1 = Z_1 \\[.25ex]
-\frac{1}{\sqrt{2}} \frac{\beta}{s}  x_6^{2} \,,
	& \mathrm{if\;} \mathfrak{n}_1 = Z_4\\
0   & \mathrm{else}
\end{cases}\\[1ex]
%
%
v_7 &= \begin{cases}
\makebox[\temp][l]{$-\frac{\alpha}{2} x_6^2 \,,$}	& \mathrm{if\;} \mathfrak{n}_1 = Z_3\\
0    & \mathrm{else}
\end{cases}\\[1ex]
\end{align*}
%
%
\textbf{b)}
Let $\mathfrak{h} = \mathbb{R}\cdot I \ltimes \mathfrak{n}_1$, where $
\mathfrak{n}_1 \in \left\lbrace \mathfrak{n},\mathfrak{n}(1,3),\mathfrak{n}(2,3),\mathfrak{n}(1,2,4),\mathfrak{n}(2,3,4),Z_1,Z_2,Z_3,Z_4,Z_5 \right\rbrace
$.
\small
\settowidth{\temp}{$x_1 x_6 \left( 2-\frac{1}{2} x_6 \right) + \alpha \left(\left(x_3 + x_5 \right) \left( x_6^{2} - x_7^{2} \right)\right) + \sqrt{2} \alpha^{2}  \left( x_6^{2} - x_7^{2} \right)$}
%
%
\begin{align*}
r_6 &= \begin{cases}
\makebox[\temp][l]{$ \sqrt{2} x_4 x_5 \cdot \mathrm{e}^{x_6 x_7} \,,$}
	& \mathrm{if\;} \mathfrak{n}_1 = \mathfrak{n}\\[.25ex]
0 \,,
	& \mathrm{if\;} \mathfrak{n}_1 = \mathfrak{n}(1,3)\\[.25ex]
 x_3 \cdot \mathrm{e}^{x_6 x_7} 	 \,,
	& \mathrm{if\;} \mathfrak{n}_1 = \mathfrak{n}(2,3)\\[.25ex]
\sqrt{2} x_4 x_5\cdot\mathrm{e}^{x_6 x_7}	\,,
	& \mathrm{if\;} \mathfrak{n}_1 = \mathfrak{n}(1,3,4)\\[.25ex]
\left(x_3 x_5 + x_4^2\right) \cdot \mathrm{e}^{x_6 x_7} \,,
	& \mathrm{if\;} \mathfrak{n}_1 = \mathfrak{n}(2,3,4)\\[.25ex]
\frac{1}{2} x_5^{2}	\,,
	& \mathrm{if\;} \mathfrak{n}_1 = Z_1 \\[.25ex]
-x_7 \left( \sqrt{2} x_4 x_6 - \left( x_3 - x_5 \right) \cdot \mathrm{e}^{\frac{1}{2} x_6 x_7^2} \right)	 \,,
	& \mathrm{if\;} \mathfrak{n}_1 = Z_2\\[.25ex]
  \left(\alpha^2 x_3 +\sqrt{2} \alpha x_4 + x_5 \right) \cdot \left( x_6 -\frac{1}{2} x_7^{2} \right) \cdot \mathrm{e}^{\frac{1}{2}x_6^2 x_7} \,,
	& \mathrm{if\;} \mathfrak{n}_1 = Z_3\\[.25ex]
-\frac{1}{2}\left(x_6^2 - 2 x_7 \right) \left(\beta x_3 - x_5 \right) \cdot \mathrm{e}^{\frac{1}{2}x_6 x_7^2} \,,
	& \mathrm{if\;} \mathfrak{n}_1 = Z_4\\[.25ex]
\frac{1}{2} x_5^2 + x_3 \cdot \left( \kappa x_6 -x_6 x_7 \right) \cdot \mathrm{e}^{\frac{1}{2}x_6 x_7^2}	\,,
	& \mathrm{if\;} \mathfrak{n}_1 = Z_5
\end{cases}\\[1ex]
%
%
r_7 &= \begin{cases}
\makebox[\temp][l]{$ x_1 x_6+  x_3^{2}+ x_5^{2} \,,$}	
	& \mathrm{if\;} \mathfrak{n}_1 = \mathfrak{n}\\[.25ex]
0 \,,
	& \mathrm{if\;} \mathfrak{n}_1 = \mathfrak{n}(1,3)\\[.25ex]
 x_1 x_6+2\sqrt{2} x_4 \,,
	& \mathrm{if\;} \mathfrak{n}_1 = \mathfrak{n}(2,3)\\[.25ex]
x_1 x_6+ x_5^{2}	\,,
	& \mathrm{if\;} \mathfrak{n}_1 = \mathfrak{n}(1,3,4)\\[.25ex]
 x_1 x_6 + 2\sqrt{2} x_4 x_5	\,,
	& \mathrm{if\;} \mathfrak{n}_1 = \mathfrak{n}(2,3,4) \\[.25ex]
\frac{1}{2} x_1 x_6^2 + \sqrt{2} x_4 x_6\cdot\mathrm{e}^{-\frac{1}{2}x_6^2 x_7}  \,,
	& \mathrm{if\;} \mathfrak{n}_1 = Z_1\\[.25ex]
 x_1 x_6 x_7 + 2\sqrt{2} x_4 x_7  - x_5 x_6 x_7 \cdot \mathrm{e}^{-\frac{1}{2} x_6 x_7^{2}}	\,,
	& \mathrm{if\;} \mathfrak{n}_1 = Z_2\\[.25ex]
\frac{1}{2} x_1 x_6^{2}  + 2 \alpha \left( x_3 + \sqrt{2} \alpha x_4 +  x_5 \right) \cdot  \left( x_6 -\frac{1}{2} x_7^{2} \right)  \,,
	& \mathrm{if\;} \mathfrak{n}_1 = Z_3\\[.25ex]
x_1 x_6 x_7 - \sqrt{2}\beta x_4 \left(x_6^2 - 2x_7\right) \,,
	& \mathrm{if\;} \mathfrak{n}_1 = Z_4\\[.25ex]
x_1 x_6 x_7 + 2\sqrt{2}\, x_4 \left(\kappa  x_6 - x_6 x_7 \right)  \,,
	& \mathrm{if\;} \mathfrak{n}_1 = Z_5
\end{cases}\\[1ex]
%
%
s_6 &= \begin{cases}
\makebox[\temp][l]{$ x_1 x_6\cdot \mathrm{e}^{x_6 x_7}+\sqrt{2} x_4 x_6 x_7 \,,$} 
	& \mathrm{if\;} \mathfrak{n}_1 = \mathfrak{n}(1,3)\\[.25ex]
0	& \mathrm{else}
\end{cases}\\[1ex]
%
%
s_7 &= \begin{cases}
\makebox[\temp][l]{$ x_2 x_6+\sqrt{2} x_3 x_4 \,,$}	
	& \mathrm{if\;} \mathfrak{n}_1 = \mathfrak{n}\\[.25ex]
 x_2 x_6+\frac{1}{2} x_3^{2}	+ x_5 x_6 x_7\cdot\mathrm{e}^{-x_6 x_7} \,,
	& \mathrm{if\;} \mathfrak{n}_1 = \mathfrak{n}(1,3)\\[.25ex]
x_2 x_6 + x_5	\,,
	& \mathrm{if\;} \mathfrak{n}_1 = \mathfrak{n}(2,3)\\[.25ex]
 x_2 x_6+\frac{1}{2} x_3^{2}	\,,
	& \mathrm{if\;} \mathfrak{n}_1 = \mathfrak{n}(1,3,4)\\[.25ex]
 x_2 x_6+\frac{1}{2} x_5^2 \,,
	& \mathrm{if\;} \mathfrak{n}_1 = \mathfrak{n}(2,3,4)\\[.25ex]
\frac{1}{2} x_2 x_6^2 + x_6 \left(  x_3+ x_5 \right) \mathrm{e}^{-\frac{1}{2} x_6^2 x_7} \,,
	& \mathrm{if\;} \mathfrak{n}_1 = Z_1 \\[.25ex]
x_2 x_6 x_7+x_5 x_7 \,,
	& \mathrm{if\;} \mathfrak{n}_1 = Z_2\\[.25ex]
\frac{1}{2} x_2 x_6^{2}  +  \left( x_3 + \sqrt{2} \alpha x_4 + \alpha^2 x_5 \right) \cdot \left( x_6 -\frac{1}{2} x_7^{2} \right) \,,
    & \mathrm{if\;} \mathfrak{n}_1 = Z_3\\[.25ex]
 x_2 x_6 x_7	+ \frac{1}{2} \left( s x_3 - \beta x_5 \right)  \left( x_6^2 - 2 x_7 \right) \,,
	& \mathrm{if\;} \mathfrak{n}_1 = Z_4\\[.25ex]
 x_2 x_6 x_7+ \left(\kappa  x_6 - x_6 x_7 \right)  \left(\frac{x_3}{\kappa} + x_5 \right) \,,
	& \mathrm{if\;} \mathfrak{n}_1 = Z_5
\end{cases}\\[1ex]
%
%
t_6 &= \begin{cases}
\makebox[\temp][l]{$\frac{1}{2} x_6 x_7^{2} \,,$}	& \mathrm{if\;} \mathfrak{n}_1 = Z_2\\[.25ex]
0
	& \mathrm{else} 
\end{cases}\\[1ex]
%
%
u_7 &= \begin{cases}
\frac{1}{\sqrt {2}} x_6^{2} x_7	\,,
	& \mathrm{if\;} \mathfrak{n}_1 = \mathfrak{n}(1,3)\\[.25ex]
\frac{1}{\sqrt {2}} x_6^{2}	\,,
	& \mathrm{if\;} \mathfrak{n}_1 = Z_1 \\[.25ex]
\makebox[\temp][l]{$0$}	& \mathrm{else} 
\end{cases}\\[1ex]
%
%
w_6 &= \begin{cases}
\makebox[\temp][l]{$ x_6 x_7 \,,$} 
	& \mathrm{if\;} \mathfrak{n}_1 \not\in \left\lbrace Z_1,\ldots,Z_5\right\rbrace \\[.25ex]
\frac{1}{2} x_6^2 x_7 \,,
	& \mathrm{if\;} \mathfrak{n}_1 \in \left\lbrace Z_1,Z_3 \right\rbrace\\[.25ex]
\frac{1}{2} x_6 x_7^2 \,,
	& \mathrm{if\;} \mathfrak{n}_1 \in \left\lbrace Z_2,Z_4,Z_5 \right\rbrace\\[.25ex]
\end{cases}\\[1ex]
%
%
t_7 &= \makebox[\temp][l]{$u_6=v_7 = 0$}\qquad\mathrm{in\; all\, cases} \\
\end{align*}